\documentclass[10pt,a4paper]{article}
\usepackage{amsfonts}
\usepackage{amssymb}
\usepackage{amsmath}
\usepackage{epsfig}
\usepackage{amsthm}
\usepackage[latin1]{inputenc}
\usepackage[english]{babel}

\newtheorem{defi}{Definition}[section]
\newtheorem{proposition}[defi]{Proposition}
\newtheorem{recall}[defi]{Recall}
\newtheorem{theoreme}[defi]{Theorem}

\newtheorem{lemme}[defi]{Lemma}
\newtheorem{corollaire}[defi]{Corollary}
\newtheorem{remarque}[defi]{Remark}

\newcommand{\F}{\mathbb F} 
\newcommand{\Z}{\mathbb Z}
\newcommand{\N}{\mathbb N}

\newcommand{\p}{\mathbb P}

\newcommand{\id}{id}

\title{Large $p$-groups actions with $\frac{|G|}{g^2} \geq \frac{4}{(p^2-1)^2}$.}
\author{Magali Rocher.}
\date{}

\begin{document}

\maketitle

 \begin{abstract}
Let $k$ be an algebraically closed field of characteristic $p>0$
and $C$ a connected nonsingular projective curve over $k$ with genus $g \geq 2$.
Let $(C,G)$ be a "big action" , i.e. a pair $(C,G)$ where $G$ is a $p$-subgroup of the $k$-automorphism group of $C$ such that$\frac{|G|}{g} >\frac{2\,p}{p-1}$. We first study finiteness results on the values taken by the quotient $\frac{|G|}{g^2}$ when $(C,G)$ runs over the big actions satisfying $\frac{|G|}{g^2} \geq M$, for a given positive real $M>0$. Then, we exhibit a classification and a parametrization of such big actions when $M=\frac{4}{(p^2-1)^2}$.
\end{abstract}

\section{Introduction.}

 \textit{Setting.} Let $k$ be an algebraically closed field of positive characteristic $p>0$ and $C$ a
connected nonsingular projective curve over $k$, with genus  $g \geq 2$. As in characteristic zero, the $k$-automorphism group of
 the curve $C$, $Aut_k(C)$, is a finite group whose order is bounded from above by a polynomial in $g$
 (cf. \cite{St} and \cite{Sin}). But, contrary to the case of characteristic zero, the bound is no more
linear but biquadratic, namely:  $|Aut_k(C)| \leq 16\, g^4$, except for the Hermitian curves: $W^q+W=X^{1+q}$, with $q=p^n$ (cf. \cite{St}). The difference  is due to the appearance of wild ramification.
More precisely, let $G$ be a subgroup of $Aut_k(C)$. If the order of $G$ is prime to $p$,
then the Hurwitz bound still holds, i.e.  $|G| \leq 84\, (g-1)$. Now, if $G$ is a $p$-Sylow subgroup of 
$Aut_k(C)$, Nakajima (cf. \cite{Na}) proves that $|G|$ can be larger according to the value of the $p$-rank $\gamma$ of the curve $C$.
Indeed, if $\gamma >0$, then $|G| \leq \frac{2\,p}{p-1} \, g$, whereas for $\gamma=0$,
 $|G| \leq \max \{g, \frac{4\,p}{(p-1)^2} \, g^2\}$ , knowing that the
quadratic upper bound $\frac{4\,p}{(p-1)^2} \, g^2$ can really be attained. 
Following Nakajima's work, Lehr and Matignon explore  the "big actions", that is to say the pairs $(C,G)$ where $G$ is a $p$-subgroup of $Aut_k(C)$ such that
$\frac{|G|}{g}> \frac{2\,p}{p-1}$ (see \cite{LM}). 
In this case, the ramification locus of the cover $\pi: C \rightarrow C/G$ is located at one point of $C$, say $\infty$. In \cite{MR}, we display necessary conditions on $G_2$, the second ramification group of $G$ at $\infty$ in lower notation, for $(C,G)$ to be a big action. In particular, we show that $G_2$ coincides with the derived subgroup $G'$ of $G$.

\medskip
\emph{Motivation and purpose.} 
The aim of this paper is to pursue the classification of big actions as initiated in \cite{LM}.
Indeed, when searching for a classification of big actions, it naturally occurs that the quotient $\frac{|G|}{g^2}$ has a "sieve" effect. Lehr and Matignon first prove that the big actions such that $\frac{|G|}{g^2} \geq \frac{4}{(p-1)^2}$ correspond to the $p$-cyclic \'etale covers of the affine line parametrized by an Artin-Schreier equation: $W^p-W=f(X):=X\,S(X)+c\,X \in k[X]$, where $S(X)$ runs over the additive polynomials of $k[X]$. In \cite{MR}, we show that the big actions satisfying $\frac{|G|}{g^2} \geq \frac{4}{(p^2-1)^2}$ correspond to the \'etale covers of the affine line with Galois group $G' \simeq (\Z/p\Z)^n$, with $n \leq 3$. This motivated the study of big actions with a $p$-elementary abelian $G'$, say $G' \simeq (\Z/p\Z)^n$, which is the main topic of \cite{MR2} where we generalize the structure theorem obtained in the $p$-cyclic case. Namely, we prove that when $G' \simeq (\Z/p\Z)^n$ with $n \geq1$, then the function field of the curve is parametrized by $n$ Artin Schreier equations: $W_i^p-W_i=f_i(X) \in k[X]$ where each function $f_i$ can be written as a linear combination over $k$ of products of at most $i+1$ additive polynomials. In this paper, we display the parametrization of the functions $f_i$'s in the case of big actions satisfying $\frac{|G|}{g^2} \geq \frac{4}{(p^2-1)^2}$. In what follows, this condition is called condition $(*)$.

\medskip
\emph{Outline ot the paper.} 
The paper falls into two main parts. The first one is focused on finiteness results for big actions $(C,G)$ satisfying $\frac{|G|}{g^2} \geq M$ for a given positive real $M>0$, called big actions satisfying $\mathcal{G}_M$,  whereas the second part is dedicated to the classification of such big actions when $M=\frac{4}{(p^2-1)^2}$.
More precisely, we prove in section 4 that, for a given $M>0$, the order of $G'$ only takes a finite number of values for $(C,G)$ a big action satisfying $\mathcal{G}_M$.  When exploring similar finiteness results for $g$ and $|G|$, we are lead to a purely group-theoretic discussion around the inclusion  $Fratt(G') \subset  [G',G]$, where $Fratt(G')$ means the Frattini subgroup  of $G'$ and $[G',G]$ denotes the commutator subgroup of $G'$ and $G$ (cf. section 4). When the inclusion is strict, $|G|$ and $g$ also take a finite number of values for $(C,G)$ satisfying $\mathcal{G}_M$. This is no more true when $Fratt(G')=[G',G]$. In this case, we can only conclude that, for $p>2$, 
the quotient $\frac{|G|}{g^2}$ takes a finite number of values for $(C,G)$ satisfying $\mathcal{G}_M$ with an abelian $G'$. Note that we do not know yet examples of big actions with a non-abelian $G'$. Another central question to is the link between the subgroups $G$ of $Aut_k(C)$ such that $(C,G)$ is a big action and a $p$-Sylow subgroup of $Aut_k(C)$ containing $G$ (section 3). Among other things, we prove that they have the same derived subgroup. This, together with the fact that the order of $G'$ takes a finite number of values for big actions satisfying $\mathcal{G}_M$, implies, on the one hand, that the order of $G'$ is a key criterion to classify big actions and, on the other hand, that we can concentrate on $p$-Sylow subgroups of $A$.
In section 5, we eventually display the classification and the parametrization of big actions $(C,G)$ under condition $(*)$ according to the order of $G'$. Pursuing the preceding discussion, we have to distinguish the cases $[G',G]=Fratt(G')(=\{e\})$ and $[G',G] \supsetneq  Fratt(G')(=\{e\}).$

\medskip
\textit{Notation and preliminary remarks.}
Let $k$ be an algebraically closed field of characteristic $p>0$.
We denote by $F$ the Frobenius endomorphism for a $k$-algebra. Then, $\wp$ means the Frobenius operator minus identity.
We denote by $k\{F\}$ the $k$-subspace of $k[X]$ generated by the polynomials $F^i(X)$, with $i \in \N$. It is a ring under the composition. Furthermore, for all $\alpha$ in $k$, $F\, \alpha=\alpha^p\,F$. The elements of $k\{F\}$ are the additive polynomials, i.e. the polynomials $P(X)$ of $k[X]$ such that for all $ \alpha$ and $\beta$ in $k$, $P(\alpha+ \beta) = P(\alpha)+ P(\beta)$.  Moreover, a separable polynomial is additive if and only if the set of its roots is a subgroup of $k$ (see \cite{Go} chap. 1).\\
\indent Let $f(X)$ be a polynomial of $k[X]$. Then, there is a unique polynomial $red (f)(X)$ in $k[X]$, called the reduced representative of $f$, which is $p$-power free, i.e. $red(f)(X) \in \bigoplus_{(i,p)=1} k\, X^i$, and such that $red(f)(X)=f(X)$ mod $\wp (k[X]).$ We say that the polynomial $f$ is reduced mod $\wp(k[X])$ if and only if it coincides with its reduced representative $red(f)$.
The equation $W^p-W=f(X)$ defines a $p$-cyclic \'etale cover of the affine line that we denote by $C_f$. 
Conversely, any $p$-cyclic \'etale cover of the affine line $Spec \, k[X]$ corresponds to a curve $C_f$ where $f$ is a polynomial of $k[X]$ (see \cite{Mi} III.4.12, p. 127). 
By Artin-Schreier theory, the covers $C_f$ and $C_{red(f)}$ define the same $p$-cyclic covers of the affine line.
The curve $C_f$ is irreducible if and only if $red(f) \neq 0$.\\
\indent Throughout the text, $C$ denotes a connected nonsingular projective curve over $k$, with genus $g \geq 2$.
We denote by $A:=Aut_kC$ the $k$-automorphism group of the curve $C$ and by $S(A)_p$ any $p$-Sylow subgroup of $A$.
 For any point $P \in C$ and any $i \geq -1$, we denote by $A_{P,i}$ the $i$-th ramification group of $A$ at $P$
in lower notation, namely
$$A_{P,i}:=\{ \sigma \in A, \, v_P( \sigma(t_{P})-t_P) \geq i+1\}$$
where $t_P$ denotes a uniformizing parameter at $P$ and $v_P$ means the order function at $P$.

\section{The setting: generalities about big actions.}

\begin{defi}
Let $C$ be a connected nonsingular projective curve over $k$, with genus $g \geq 2$.
Let $G$ be a subgroup of $A$.
We say that the pair $(C,G)$ is a big action if $G$ is a finite $p$-group such that
$$
\frac{ |G|}{g} > \frac{2\,p}{p-1}
$$
\end{defi}

\noindent To precise the background of this work, we first recall basic properties of big actions established in \cite{LM} and \cite{MR}. 

\begin{recall} 
Assume that $(C,G)$ is a big action. Then, there is a point of $C$ (say $\infty$) such that $G$ is the wild inertia subgroup of $G$ at $\infty$: $G_1$.
Moreover, the quotient $C / G$ is isomorphic to the projective line $\p^1_k$ and the ramification locus (respectively branch locus) of the cover $\pi: \, C \rightarrow C/G$ is the 
point  $\infty$ (respectively $\pi(\infty))$.
For all $i \geq 0$, we denote by $G_i$ the $i$-th lower ramification group of $G$ at $\infty$:
$$G_i:=\{ \sigma \in G, \, v_{\infty}( \sigma(t_{\infty})-t_{\infty}) \geq i+1\}$$
where $t_{\infty}$ denotes a uniformizing parameter at $\infty$ and $v_{\infty}$ means the order function at $\infty$.
\begin{enumerate}
\item Then, $G_2$ is non trivial and it is strictly included in $G_1$.
\item The quotient curve $C/G_2$ is isomorphic to the projective line $\p_k^1.$ 
\item The quotient group $G /G_2$ acts as a group of translations of the affine line $C/G_2 - \{\infty \}=Spec \,k[X]$, through $X  \rightarrow X+y$, where $y$ runs over a subgroup $V$ of $k$. Then, $V$ is an
$\F_p$-subvector space of $k$. We denote by $v$ its dimension. This gives the following exact
  sequence:
   $$ 0 \longrightarrow G_2 \longrightarrow G=G_1  \stackrel{\pi}{\longrightarrow} V \simeq
   (\Z/\,p\, \Z)^v \longrightarrow 0$$ where
 $$  \pi: 
\left\{
\begin{array}{ll}
 G \rightarrow V\\
g \rightarrow g(X)-X
\end{array}
\right.
$$ 
\end{enumerate}
\end{recall}

\begin{recall} \label{rappel} (\cite{MR} Thm. 2.6.4).
Let $(C,G)$ be a big action. Then, $$G_2=G'=Fratt(G)$$ where $G'$ means the commutator subgroup of $G$ and $Fratt(G)=G'G^p$ the Frattini subgroup of $G$.
\end{recall}

\noindent To conclude this first section, we introduce new definitions used in our future classification.

\begin{defi} \label{def11}
Let $C$ be a connected nonsingular projective curve over $k$, with genus $g \geq 2$.
Let $G$ be a subgroup of $A$.
Let $M  > 0$ be a positive real. 
We say that:
\begin{enumerate}
\item $G$ satisfies $\mathcal{G}(C)$ (or $(C,G)$ satisfies $\mathcal{G}$) if $(C,G)$ is a big action.
\item $G$ satisfies $\mathcal{G}_M(C)$ (or $(C,G)$ satisfies $\mathcal{G}_M$) if $(C,G)$ is a big action with $\frac{|G|}{g^2} \geq M.$ 
\item If $(C,G)$ satisfies $\mathcal{G}_M$ with $M=\frac{4}{(p^2-1)^2}$, we say that $(C,G)$ satisfies condition $(*)$.
\end{enumerate}
\end{defi}

\begin{remarque} \label{rem10}
There exists big actions $(C,G)$ satisfying $\mathcal{G}_M$ if and only if $M \leq \frac{4\,p}{(p-1)^2}$ (see \cite{St}).
\end{remarque}

\section{A study on $p$-Sylow subgroups of $Aut_k(C)$ inducing big actions.}

\noindent In this section, we more specifically concentrate on the $p$-Sylow subgroup(s) of $A$ satisfying $\mathcal{G}(C)$ (resp. $\mathcal{G}_M(C)$).

\begin{remarque} \label{italien} 
Let $C$ be a connected nonsingular projective curve over $k$, with genus $g \geq 2$.
Assume that there exists a subgroup $G \subset A$ satisfying $\mathcal{G}(C)$.
\begin{enumerate}
\item Then, every $p$-Sylow subgroup of $A$ satisfies $\mathcal{G}(C)$.
\item Moreover, $A$ has a unique $p$-Sylow subgroup except in the three following cases (cf. \cite{Han} and \cite{Gi}):
\begin{enumerate}
\item The Hermitian curve $$C_H : \quad  W^q+W=X^{1+q}$$ with $p\geq 2$, $q=p^s$, $s \geq 1$.
Then, $g=\frac{1}{2} \,(q^2-q)$
and $A \simeq PSU(3,q)$ or $A \simeq PGU(3,q)$. It follows that $|A|=q^3\,(q^2-1)\,(q^3+1)$, so $\frac{|S(A)_p|}{g} =\frac{2q^2}{q-1} >\frac{2\,p}{p-1}$ and $\frac{|S(A)_p|}{g^2} =\frac{4\,q}{(q-1)^2}$, where $S(A)_p$ denotes any $p$-Sylow subgroup of $A$. Thus, $(C_H,S(A)_p)$ is a big action with $G'=G_2 \simeq (\Z/p\Z)^s$. It satisfies condition $(*)$ if and only if $1 \leq s \leq 3$.\\

\item The Deligne-Lusztig curve arising from the Suzuki group $$C_S: \quad W^q+W=X^{q_0} \, (X^q+X)$$ with $p=2$, 
$q_0=2^s$, $s \geq 1$ and $q=2^{2s+1}$.
 In this case, $g= q_0(q-1)$ and $A \simeq Sz(q)$ the Suzuki group.
It follows that $|A|=q^2\,(q-1)\,(q^2+1)$, so $\frac{|S(A)_p|}{g} =\frac{q^2}{q_0(q-1)}>\frac{2\,p}{p-1}$ and $\frac{|S(A)_p|}{g^2} =\frac{ q^2}{q_0^2\,(q-1)^2} < \frac{4}{(p^2-1)^2}$, for all $s \geq 1$. 
Thus, $(C_S,S(A)_p)$ is a big action with $G'=G_2 \simeq (\Z/p\Z)^{2s+1}$ but it never satisfies condition $(*)$.\\

\item The Deligne-Lusztig curve arising from the Ree group
$$C_R: \quad W_1^q-W_1=X^{q_0} \, (X^q+X) \quad \mbox{and} \quad W_2^q-W_2=X^{2q_0} \, (X^q+X)$$ with $p=3$, 
$q_0=3^s$, $s \geq 1$ and $q=3^{2s+1}$. 
Then, $g= \frac{3}{2} \,q_0\,(q-1)\, (q+q_0+1)$ and $A \simeq Ree(q)$ the Ree group.
It follows that $|A|=q^3\,(q-1)\,(q^3+1)$, so
$\frac{|S(A)_p|}{g} =\frac{2q^3}{3q_0(q-1)(q+q_0+1)}>\frac{2\,p}{p-1}$ and 
$\frac{|S(A)_p|}{g^2} =\frac{4q^3}{9 q_0^2(q-1)^2(q+q_0+1)^2}< \frac{4}{(p^2-1)^2}$ for all $s \geq 1$.
 Thus, $(C_R,S(A)_p))$ is a big action with $G'=G_2 \simeq (\Z/p\Z)^{2(2s+1)}$ but it never satisfies condition $(*)$.
\end{enumerate}
In each of these three cases, the group $A$ is simple, so it has more than one $p$-Sylow subgroups.
\end{enumerate}
\end{remarque}

\noindent Now, fix $C$ a connected nonsingular projective curve over $k$, with genus $g \geq 2$. We  highlight the link between the groups $G$ satisfying $\mathcal{G}(C)$ (resp. $\mathcal{G}_M(C)$) and the $p$-Sylow subgroup(s) of $A$.

\begin{proposition} \label{prop11}
Let $C$ be a connected nonsingular projective curve over $k$, with genus $g \geq 2$.
\begin{enumerate}
\item Let $G$ satisfy $\mathcal{G}(C)$.  
\begin{enumerate}
\item Then, there exists a point of $C$, say $\infty$, such that $G$ is included in
$A_{\infty,1}$. For all $i \geq 0$, we denote by $G_i$ the $i$-th ramification group of $G$ at $\infty$ in lower notation. Then, $A_{\infty,1}$ satisfies $\mathcal{G}(C)$ and $A_{\infty,2}=G_2$, i.e. $(A_{\infty,1})'=G'$.
Thus, we obtain the following diagram:
$$
\begin{array}{ccccccccc}
0 &\longrightarrow &A_{\infty,2} &
\longrightarrow &A_{\infty,1}&\stackrel{\pi}{\longrightarrow}& W \subset k& \longrightarrow &0 \\
& &|| & & \cup & & \cup  & &\\
0 &\longrightarrow &G_2 &
\longrightarrow &G=G_1&\stackrel{\pi}{\longrightarrow}& V& \longrightarrow &0 
\end{array}
$$
In particular, $G=\pi^{-1}(V)$ where $V$ is an $\F_p$-subvector space of $W$.
\item $A_{\infty,1}$ is a $p$-Sylow subgroup of $A$. Moreover, except in the three special cases mentionned in Remark \ref{italien}, $A_{\infty,1}$ is  the unique $p$-Sylow subgroup of $A$.
\item  Let $M$ be a positive real such that $G$ satisfies $\mathcal{G}_M(C)$. Then, $A_{\infty,1}$ also satifies $\mathcal{G}_M(C).$
\end{enumerate}
\item Conversely, let $\infty$ be a point of the curve $C$ such that $A_{\infty,1}$ satisfies $\mathcal{G}(C)$.
Consider $V$ an $\F_p$-vector space of $W$, defined as above, and put $G:=\pi^{-1}(V)$. 
\begin{enumerate}
\item Then, the group $G$ satisfies $\mathcal{G}(C)$ if and only if
$$|W| \geq |V| > \frac{2\,p}{p-1} \, \frac{g}{|A_{\infty,2}|}$$
\item  Let $M$ be a positive real such that $A_{\infty,1}$ satisfies $\mathcal{G}_M(C)$. Then, $G$ satisfies  $\mathcal{G}_M(C)$ if and only if 
$$ |W| \geq |V|  \geq M \, \frac{g^2}{|A_{\infty,2}|}$$ 
\end{enumerate}
\end{enumerate}
\end{proposition}

\textbf{Proof:} 
The first assertion (1.a) derives from \cite{LM} (Prop 8.5) and \cite{MR} (Cor. 2.10). The second point (1.b) comes from \cite{MR} (Rem 2.11) together with Remark \ref{italien}. The other claims are obtained via calculation. $\square$

\begin{remarque}
Except in the three special cases mentionned in Remark \ref{italien}, the point $\infty$ of $C$ defined in Proposition 
\ref{prop11} is uniquely determined. In particular, except for the three special cases, if $P$ is a point of $C$ such that $A_{P,1}$ satisfies $\mathcal{G}(C)$, then $P=\infty$.
\end{remarque}

As a conclusion, if $G$ satisfyies $\mathcal{G}(C)$ (resp. $\mathcal{G}_M(C)$)  and if $A_{\infty,1}$ is a (actually "the", in most cases) $p$-Sylow subgroup of $A$ containing $G$, then $A_{\infty,1}$ also satisfies $\mathcal{G}(C)$ (resp. $\mathcal{G}_M(C)$) and has the same derived subgroup. So, in our attempt to classify the big actions $(C,G)$ satisfying $\mathcal{G}_M$, this leads us to focus on the derived subgroup $G'$ of $G$. 

\section{Finiteness results for big actions satisfying $\mathcal{G}_M$.}

\subsection{An upper bound on $|G'|$.}

\begin{lemme} \label{finite}
Let $M>0$ be a positive real such that $(C,G)$ is a big action satisfying $\mathcal{G}_M$. Then, the order of $G'$ is bounded as follows:
$$p \leq |G'| \leq \,\frac{4\,p}{(p-1)^2} \,\frac{2+M+2\,\sqrt{1+M}}{M^2}$$ Thus, $|G'|$ only takes a finite number of values for $(C,G)$ a big action satisfying $\mathcal{G}_M$.
\end{lemme}

\noindent \textbf{Proof:}
We first recall that $G'=G_2$ is a non-trivial $p$-group (see e.g. \cite{LM} Prop. 8.5).
Now, let $i_0 \geq 2$ be the integer such that the lower ramification filtration of $G$ at $\infty$ reads:
$$G=G_0=G_1 \supsetneq G_2= \cdots =G_{i_0} \supsetneq G_{i_0+1} = \cdots$$
Put $|G_2/G_{i_0+1}|=p^m$, with $m \geq 1$, and
 $\mathcal{B}_m := \frac{4}{M}  \frac{|G_2/ G_{i_0+1}|}{(|G_2/ G_{i_0+1}|-1)^2}=\frac{4}{M} \frac{p^{m}}{(p^{m}-1)^2}$. By \cite{LM} (Thm. 8.6),  $M \leq \frac{|G|}{g^2}$ implies $1<|G_2| \leq  \frac{4}{M}  \frac{|G_2/ G_{i_0+1}|^2}{( |G_2/ G_{i_0+1}|-1)^2} =p^m \, \mathcal{B}_m$. 
 From $|G_2| =p^{m} |G_{i_0+1}|$, we infer $1 \leq |G_{i_0+1}| \leq \mathcal{B}_m$.
Since $(\mathcal{B}_m)_{m \geq 1}$ is a decreasing sequence which tends to $0$ as $m$ grows large, we conclude that $m$ is bounded. More precisely, $m < m_0$ where $m_0$ is the smallest integer such that $\mathcal{B}_{m_0} <1$.  As $M \leq \frac{4\,p}{(p-1)^2} \leq 8$ (see Remark \ref{rem10}), computation shows that
$\mathcal{B}_m<1  \Leftrightarrow  p^m > \phi(M):=\frac{2+M+2\,\sqrt{1+M}}{M}$.
As $(\mathcal{B}_m)_{m \geq 1}$ is decreasing, 
$$|G_2| \leq p^m \, \mathcal{B}_m \leq \phi_1(M) \, \mathcal{B}_1= \frac{\phi(M)}{M} \, \frac{4\,p}{(p-1)^2}$$ 
The claim follows. $\square$

\medskip

\noindent We deduce that, for big actions $(C,G)$ satisfying $\mathcal{G}_M$, an upper bound on $|V|$ induces an upper bound on the genus $g$ of $C$. 

\begin{corollaire} \label{borneV}
Let $M>0$ be a positive real such that $(C,G)$ is a big action satisfying $\mathcal{G}_M$. Then, 
$$g < \frac{|G'|\, |V|}{2\,p} \leq  \frac{2}{p-1}\, \frac{2+M+2\sqrt{1+M}}{M^2}\, |V|$$
\end{corollaire}

This raises the following question.
Let $(C,G)$ be a big action satisfying $\mathcal{G}_M$; in which cases is $|V|$ (and then $g$) bounded from above? In other words, in which cases, does the quotient $\frac{|G|}{g}$ take a finite number of values when $(C,G)$ satisfy $\mathcal{G}_M$?
We begin with preliminary results on big actions leading to a purely group-theoretic discussion leading to compare the Frattini subgroup of $G'$ with the commutator subgroup of $G'$ and $G$.

\subsection{Preliminaries to a group-theoretic discussion.}

\begin{lemme} \label{pgroup}
Let $(C,G)$ be a big action. If $G' \subset Z(G)$, then
$G'(=G_2)$ is $p$-elementary abelian, say $G'\simeq (\Z/p\Z)^n$, with $n \geq 1$.
In this case, the function field $L=k(C)$ is parametrized by $n$ equations:
$$\forall \, i \in \{1,\cdots,n\}, \quad W_i^p-W_i=f_i(X)=X\, S_i(X)+c_i\, X \in k[X]$$
where $S_i$ is an additive polynomial of $k[X]$ with degree $s_i \geq 1$ in $F$ and $s_1 \leq s_2\cdots \leq s_n.$
Moreover, $V \subset \cap_{1\leq i \leq n} Z(Ad_{f_i})$ where $Ad_{f_i}$ denotes the palindromic polynomial related to $f_i$ as defined in \cite{MR2} (Prop. 2.13)
\end{lemme}

\noindent \textbf{Proof:}
The hypothesis first requires $G'=G_2$ to be abelian. Now, assume that $G_2$ has exponent strictly greater than $p$. Then, there exists a surjective map $\phi: G_2 \rightarrow \Z/p^2\Z$. So $H:=Ker \phi \subsetneq G_2 \subset Z(G)$  is a normal subgroup of $G$. It follows from \cite{MR} (Lemma 2.4) that the pair $(C/H,G/H)$ is a big action with second ramification group $(G/H)_2\simeq \Z/p^2\Z$. This contradicts \cite{MR} (Thm. 5.1). The last part of the lemma comes from \cite{MR2} (Prop. 2.13).
$\square$

\begin{corollaire} \label{Fratt}
Let $(C,G)$ be a big action. Let $H:=[G',G]$ be the commutator subgroup of $G'$ and $G$.
\begin{enumerate}
\item Then, $H$ is trivial if and only if $G'\subset Z(G)$.
\item The group $H$ is strictly included in $G'$.
\item The pair $(C/H,G/H)$ is a big action. Moreover, its second ramification group
$(G/H)_2=(G/H)'=G_2/H \subset Z(G/H)$ is $p$-elementary abelian.
\end{enumerate}
\end{corollaire}

\noindent \textbf{Proof:} \begin{enumerate}
\item The first assertion is clear.
\item  As $G'$ is normal in $G$, then $H \subset G'$. Assume that $G'=H$. Then, the lower central series of $G$ is stationnary, which contradicts the fact that the $p$-group $G$ is nilpotent (see e.g. \cite{Su2} Chap.4). So $H \subsetneq G'$. 
\item As $H \subsetneq G'=G_2$ is normal in $G$, it follows from \cite{MR} (Lemma 2.4 and Thm. 2.6) that the pair $(C/H,G/H)$ is a big action with second ramification group  $(G/H)_2=G_2/H$. From $H=[G_2:G]$, we gather that $G_2/H \subset Z(G/H)$. Therefore, we deduce from Lemma \ref{pgroup} that $(G/H)_2$ is $p$-elementary abelian. $\square$
\end{enumerate}

\begin{corollaire} \label{Frattbis}
Let $(C,G)$ be a big action. Let $F:=Fratt(G')$ be the Frattini subgroup of $G'$.
\begin{enumerate}
\item Then, $F$ is trivial if and only if $G'$ is an elementary abelian $p$-group.
\item We have the following inclusions: $F \subset [G',G] \subsetneq G'$.
\item The pair $(C/F,G/F)$ is a big action. Moreover, its second ramification group
$(G/F)_2=(G/F)'=G_2/F$ is $p$-elementary abelian.
\item Let $M$ be a positive real. If $(C,G)$ satisfies $\mathcal{G}_M$, then $(C/F,G/F)$ also satisfies $\mathcal{G}_M$.
\end{enumerate}
\end{corollaire}

\noindent \textbf{Proof:} \begin{enumerate}
\item As $G'$ is a $p$-group,  $F=(G')'(G')^p$, where $(G')'$ means the derived subgroup of $G'$ and $(G')^p$ the subgroup generated by the $p$ powers of elements of $G'$ (cf. \cite{LG} Prop. 1.2.4). This proves that if $G'$ is $p$-elementary abelian, then $F$ is trivial. The converse derives from the fact that $G'/F$ is $p$-elementary abelian (cf. \cite{LG} Prop. 1.2.4).
\item Using Corollary \ref{Fratt}, the only inclusion that remains to show is $F \subset [G',G]$.  As $G'/[G',G]$ is abelian , $(G')' \subset [G',G]$. As $G'/[G',G]$ has exponent $p$, $(G')^p \subset [G',G]$. The claim follows.
\item Since $F\subsetneq G'=G_2$ is normal in $G$, we deduce from \cite{MR} (Lemma 2.4) that 
the pair $(C/F,G/F)$ is a big action with second ramification group: $(G/ F)_2=G_2/F=(G/F)'$. 
Furthermore, as $G_2$ is a $p$-group, $G_2/F$ is an elementary abelian $p$-group (see above).
\item This derives from \cite{LM} (Prop. 8.5 (ii)). $\square$
\end{enumerate}

\medskip

\noindent This leads us to discuss according to whether $Fratt(G') \subsetneq [G',G]$ or $Fratt(G') = [G',G]$.

\subsection{Case: $Fratt(G') \subsetneq [G',G]$}
\noindent We start with the special case $\{e\}=Fratt(G') \subsetneq [G',G]$,  i.e. $G'$ is $p$-elementary abelian and $G' \not \subset Z(G)$.

\begin{proposition} \label{autreborne}
Let $M>0$ be a positive real such that $(C,G)$ is a big action satisfying $\mathcal{G}_M$.
Suppose that $\{e\}=Fratt(G') \subsetneq [G',G]$. Then, $|V|$ and $g$ are bounded as follows:
\begin{equation} \label{VVV}
|V| \leq \frac{4}{M} \frac{|G_2|}{(p-1)^2} \leq \frac{16\,p}{(p-1)^4} \,\frac{2+M+2\,\sqrt{1+M}}{M^3} 
\end{equation}
and
\begin{equation} \label{ggg}
\frac{p-1}{2} \, |V| \leq g < \frac{32\,p}{(p-1)^5} \, \frac{(2+M+2\sqrt{1+M})^2}{M^5}
\end{equation}
Thus, under these conditions, $g$, $|V|$ and so the quotient $\frac{|G|}{g}$ only take a finite number of values.
\end{proposition}

\noindent \textbf{Proof:}
Write $G'=G_2 \simeq (\Z/p\Z)^n$, with $n \geq 1$.
As $G_2 \not \subset Z(G)$, \cite{MR2} (Prop. 2.13) ensures the existence of a
 smaller integer $j_0 \geq 1$ such that $f_{j_0+1}(X)$ cannot be written as $c\,X+XS(X)$, with $S$ in $k\{F\}$. If $j_0 \geq 2$, it follows that, for all $y$ in $V$, the coefficients of the matrix
$L(y)$ satisfy $\ell_{j,i}(y)=0$ for all $2\leq i \leq j_0$ and $1 \leq j \leq i-1$. Moreover, the matricial
multiplication proves that, for all $i$ in $\{1,\cdots, j_0\}$, the functions $\ell_{i,j_0+1}$ are nonzero linear forms from  $V$ to $\F_p$. 
Put $\mathcal{W}:= \bigcap_{1 \leq i \leq j_0} ker \,
 \ell_{i,j_0+1}$. 
 Let $C_{f_{j_0+1}}$ be the curve parametrized by $W^p-W=f_{{j_0}+1}(X)$. It defines an \'etale cover of the
affine line with group $\Gamma_0 \simeq \Z/p\Z$. Since, for all $y$ in $\mathcal{W}$, $f_{i_{0}+1}(X+y)= f_{i_{0}+1}(X)$ mod $\wp(k[X])$, the group of translations of the affine
   line: $\{ X \rightarrow X+y$, $y \in \mathcal{W}\} $ can be extended to a $p$-group of automorphisms of
the curve $C_{f_{j_0+1}}$ , say $\Gamma$, with the following exact sequence: 
$$ 0 \longrightarrow  \Gamma_0 \simeq \Z / p\, \Z  \longrightarrow \Gamma \longrightarrow  \mathcal{W} \longrightarrow  0 $$
 The pair $(C_{f_{j_0+1}},\Gamma)$ is not a big action. Otherwise, its second ramification group would be $p$-cyclic, which contradicts the form of the function $f_{j_0+1}(X)$, as compared with \cite{MR} (Prop. 2.5). 
Thus, $\frac{|\Gamma|}{g_{C_{f_{j_0+1}}}}=\frac{2 \, p}{p-1} \,
\frac{|\mathcal{W}|}{(m_{j_0+1}-1)} \leq \frac{2\, p}{p-1}$. 
The inequality $ \frac{|V|}{p^{j_0}} \leq |\mathcal{W}| \leq
   (m_{j_0+1}-1)$ combined with the formula given in \cite{MR2} (Cor. 2.7) yields a lower bound on the genus, namely:
   $$g=\frac{p-1}{2}\,\sum_{i=1}^n\, p^{i-1} \,(m_i-1)
   \geq \frac{p-1}{2} \, p^{j_0} \, (m_{j_0+1}-1) \geq \frac{p-1}{2} \,|V|.$$
   It follows that $M \leq \frac{|G|}{g^2}=\frac{|G_2| \, |V|}{g^2}\leq \frac{4\,|G_2|}{(p-1)^2\, |V|}$.
Using Lemma \ref{finite}, we gather inequality \eqref{VVV}. Inequality 
\eqref{ggg} then derives from Corollary \ref{borneV}. $\square$

\medskip

\noindent The following corollary generalizes the finiteness result of Proposition \ref{autreborne} to all big actions
 satisfying $\mathcal{G}_M$ such that $Fratt(G') \subsetneq [G',G]$.

\begin{corollaire}
Let $M>0$ be a positive real such that $(C,G)$ is a big action satisfying $\mathcal{G}_M$. Suppose that $Fratt(G') \subsetneq [G',G]$. Then, $|V|$ and $g$ are bounded as in Proposition \ref{autreborne}. So the quotients $\frac{|G|}{g}$ and $\frac{|G|}{g^2}$ only take a finite number of values.
\end{corollaire}

\noindent \textbf{Proof:} Put $F:=Fratt(G')$.
Corollary \ref{Frattbis} asserts that 
the pair $(C/F,G/F)$ is a big action satisfying $\mathcal{G}_M$ whose second ramification group: $(G/ F)_2=G_2/F$ is $p$-elementary abelian. From $F \subsetneq [G_2:G]$, we gather  $\{e\} \subsetneq [G_2/F: G/F]$, which implies  $(G/F)_2=(G/F)' \not \subset Z(G/F)$. We deduce that $|V|$ is bounded from above as in Proposition \ref{autreborne}. The claim follows. $\square$

\subsection{Case: $Fratt(G') = [G',G]$}

It remains to investigate the case where $Fratt(G')=[G',G]$. In particular, this equality is satisfied when $G'$ is included in the center of $G$ and so is $p$-elementary abelian (cf. Lemma 3.3), i.e. 
 $\{e\}=Fratt(G')=[G',G]$. The finiteness result on $g$ obtained in the preceding section is no more true in this case, as
illustrated by the remark below.

\begin{remarque} \label{pasdeborne}
For any integer $s\geq 1$, Proposition 2.5 in \cite{MR} exhibits an example of big actions $(C,G)$ with $C: W^p-W=X\, S(X)$ where $S$ is
an additive polynomial of $k[X]$ with degree $p^s$. In this case, $g=\frac{p-1}{2} \, p^s$, $V=Z(Ad_f)\simeq (\Z/p\Z)^{2s}$ and $G'=G_2\simeq \Z/p\Z \subset Z(G)$. It follows that $\frac{|G|}{g^2}=\frac{4\,p}{(p-1)^2}$.
So, for all $M \leq \frac{4\,p}{(p-1)^2}$, $(C,G)$ satisfies $\mathcal{G}_M$, with $\{e\}=Fratt(G')=[G',G]$, whereas $g=\frac{p-1}{2} \, p^s$ grows arbitrary large with $s$.
\end{remarque}

 Therefore, in this case, neither $g$ nor $|V|$ are bounded. Nevertheless, the following section shows that, under these conditions, the quotient $\frac{|G|}{g^2}$ take a finite number of values.

\subsubsection{Case: $Fratt(G') = [G',G]=\{e\}$.}

\begin{proposition}
Let $M>0$ be a positive real such that $(C,G)$ is a big action satisfying $\mathcal{G}_M$. Assume that $[G',G]=Fratt(G')=\{e\}$. Let $s_1$ be the integer in Lemma \ref{pgroup}. Then, $\frac{p^{2s_1}}{g^2}$ and $\frac{|V|}{p^{2s_1}}$ are bounded as follows: 
\begin{equation} \label{g22}
 \frac{p^{2s_1}}{g^2} \leq \frac{(p-1)^2}{4\,p} \, \frac{M^3}{2+M+2\, \sqrt{1+M}}
\end{equation}
and
\begin{equation} \label{V22}
 1 \leq \frac{|V|}{p^{2s_1}} \leq \frac{(p-1)^4}{16\,p} \, \frac{M^3}{2+M+2\, \sqrt{1+M}} 
\end{equation}
Thus, the quotient $\frac{|G|}{g^2}$ takes a finite number of values.
\end{proposition}

\noindent \textbf{Proof:} Write $G'=G_2\simeq (\Z/p\Z)^n$, with $n\geq 1$.
Lemma \ref{finite} first implies that $p^n$ can
 only take a finite number of values. Moreover, as recalled in Lemma \ref{pgroup}, $V \subset \bigcap_{i=1}^n Z (Ad_{f_i})$ and $|G|=|G_2||V| \leq
 p^{n+2\,s_1}$. We compute the genus by means of \cite{MR2} (Cor. 2.7):
$$g=\frac{p-1}{2} \,
 \sum_{i=1}^n p^{i-1}\, (m_i-1)=\frac{p-1}{2}\, p^{s_1}\,
  (\sum_{i=1}^n p^{i-1} p^{s_i-s_1})$$ It follows that:
 $0 <M \leq \frac{|G|}{g^2} \leq \frac{4\,p^n}{(p-1)^2 (\sum_{i=1}^n p^{i-1} p^{s_i-s_1})^2}$
 This implies $(\sum_{i=1}^n p^{i-1} p^{s_i-s_1})^2 \leq \frac{4\, p^n}{M\,
 (p-1)^2}$. As $p^n$ is bounded from above, the set $\{s_i-s_1, i \in [1,n] \} \subset \N$ is also bounded, and then
 finite.
More precisely, we gather that
$$\frac{g^2}{p^{2s_1}} = \frac{(p-1)^2}{4} \, (\sum_{i=1}^n p^{i-1} p^{s_i-s_1})^2) \leq \frac{p^n}{M}$$
Combined with Lemma \ref{finite}, this gives inequality \eqref{g22}. 
 Besides, from $ M \leq \frac{|G|}{g^2}= \frac{|V|\, p^n}{g^2}$, we infer that $\frac{1}{|V|} \leq
 \frac{p^n}{M\, g^2}$, which involves: $$1 \leq \frac{p^{2s_1}}{|V|} \leq \frac{p^{2\,s_1}\, p^n}{M\,  g^2}=
 \frac{4 \, p^n}{M\, (p-1)^2\,  (\sum_{i=1}^n p^{i-1} p^{s_i-s_1})^2} \leq \frac{4\,p^n}{M\, (p-1)^2} $$
This, together with Lemma \ref{finite}, yields inequality \eqref{V22}. In particular, the set $\{ \frac{p^{2s_1}}{|V|}  \} \subset \N$ is bounded, and then finite, as well as the set $\{ \frac{|V|}{p^{2s_1}}  \}.$ 
Therefore, the quotient $\frac{|G|}{g^2}=p^n \, \frac{|V|}{p^{2s_1}} \, \frac{p^{2s_1}}{g^2}$ can only take a finite number of values.
 $\square$

\medskip

\noindent The last remaining case is $Fratt(G')=[G',G] \neq \{e\}$.

\subsubsection{Case: $Fratt(G') = [G',G]\neq \{e\}$.}
As shown below, this case can only occur for $G'(=G_2)$ non abelian. Note that we do not know yet examples of big actions with a non abelian $G'(=G_2)$.

\begin{theoreme} \label{theo}
Assume that $p>2$. Let $(C,G)$ be a big action with $Fratt(G') = [G',G]\neq \{e\}$. Then, $G'(=G_2)$ is non abelian.
\end{theoreme}

\noindent We deduce the following 

\begin{corollaire}
Assume that $p>2$. Let $M>0$ be a positive real. Let $(C,G)$ be a big action satisfying $\mathcal{G}_M$ with $G'$ abelian. Then, $\frac{|G|}{g^2}$ only takes a finite number of values.
\end{corollaire}

\begin{remarque}
Theorem \ref{theo} is no more true for $p=2$. A counterexample is given by \cite{MR} (Prop. 6.9)
applied with $p=2$.
Indeed, when keeping the notations of \cite{MR} (Prop. 6.9), take $q=p^{e}$ with $p=2$, $e=2s-1$ and $s \geq 2$. Put $K=\F_{q}(X)$.
Let $L:=\F_{q}(X,W_1,V_1,W_2)$ be the extension of $K$ parametrized by
$$W_1^{2^{2s-1}}-W_1=X^{2^{s-1}}\, (X^{2^{2s-1}}-X) \quad \quad V_1^{2^{2s-1}}-V_1=X^{2^{s-2}}\, (X^{2^{2s-1}}-X)$$
$$[W_1,W_2]^2-[W_1,W_2]=[X^{1+2^s},0]-[X^{1+2^{s-1}},0]$$
Let $G$ be the $p$-group of $\F_{q}$-automorphisms of $L$ constructed as in \cite{MR} (Prop. 6.9.3).
Then, the formula established for $g_L$ in \cite{MR} (Prop. 6.9.4) shows that
the pair $(C,G)$ is a big action as soon as $s \geq 4$. In this case, $G'=G_2 \simeq \Z/2^2\Z \times (\Z/2\Z)^{6s-4}$ (cf. \cite{MR} Prop. 6.7.2). As the functions $X^{2^{s-1}}\, (X^{2^{2s-1}}-X)$ and $X^{2^{s-2}}\, (X^{2^{2s-1}}-X)$ are products of two additive polynomials, it follows from next proof (cf. point 6)  that $[G',G]=Fratt(G') \neq \{e\}$.
\end{remarque}

\noindent \textbf{Proof of Theorem \ref{theo}:} \begin{enumerate}
\item \emph{Preliminary remarks: the link with Theorem 5.1  in \cite{MR}.} 
\begin{enumerate}
\item One first remarks that Theorem \ref{theo} implies Theorem 5.1 in \cite{MR}. The latter states that there is no big action $(C,G)$ with $G_2$ cyclic of exponent strictly greater than $p$. Indeed, assume that there exists one.
Then, $G'=G_2$ is abelian and $Fratt(G')=(G')^p \neq \{e\}$. To contradicts Theorem \ref{theo}, it remains to show that $F:=Fratt(G')=[G',G]$.
From Corollary \ref{Frattbis}, we infer that $(C/F, G/F)$ is a big action whose second ramification group $G_2/F$ is cyclic of order $p$. Then, $(G/F)'=(G/F)_2=G_2/F \subset Z(G/F)$ (cf. \cite{MR} Prop. 2.5 and \cite{MR2} Prop. 2.13). It follows that $Fratt((G/F)')=[(G/F)', G/F]=\{e\}$. As $F \subset G'$, this imposes $F=[G',G]$. Then, Theorem \ref{theo} contradicts the fact that $G'=G_2$ is abelian.

\item The object of Theorem \ref{theo} is to prove that there exists no big action $(C,G)$ with $G'=G_2$ abelian of exponent strictly greater than $p$ such that $Fratt(G')=[G',G]$. 
The proof follows the same canvas as the proof of \cite{MR} (Thm. 5.1). Nevertheless, to refine the arguments, we use the formalism related to the ring filtration of $k[X]$ linked with the additive polynomials as introduced in \cite{MR2} (cf. section  3). More precisely, we recall that, for any $t \geq 1$, we define $\Sigma_t$ as the $k$-subvector space of $k[X]$ generated by $1$ and the products of at most $t$ additive polynomials of $k[X]$ (cf. \cite{MR2} Def. 3.1). In what follows, we assume that there exists a big action $(C,G)$ with $G'=G_2$ abelian of exponent strictly greater than $p$ such that $Fratt(G')=[G',G]$. 
\end{enumerate}

\item \emph{One can suppose that $G'=G_2 \simeq \Z/p^2\Z \times (\Z/p\Z)^r$, with $r \geq 1$.} \\
Indeed, write $G'/(G')^{p^2} \simeq (\Z/p^2\Z)^a \times (\Z/p\Z)^b$. By assumption, $a\geq 1$. Using \cite{Su1} (Chap.2, Thm. 19), one can find an index $p$-subgroup of $(G')^p$, normal in $G$, such that 
$(G')^{p^2} \subset H \subsetneq (G')^p \subsetneq G'=G_2$. Then, we infer from \cite{MR} (Lemma 2.4) that $(C/H, G/H)$ is a big action with second ramification group $(G/H)'=(G/H)_2=G_2/H \simeq (\Z/p^2\Z) \times (\Z/p\Z)^{a+b-1}$. Furthermore, as $G'$ is abelian, $Fratt(G')=(G')^p$ (resp. $Fratt((G/H)')=((G/H)')^p$). From $H \subset (G')^p$ with $H$ normal in $G$ and $Fratt(G')=[G',G]$, we gather that $Fratt((G/H)')=(G')^p/H= Fratt(G')/H=[(G/H)',G/H]$.

\item \emph{Notation.} \\
In what follows, we denote by $L:=k(C)$ the function field of $C$ and by $k(X):=L^{G_2}$ the subfield of $L$ fixed by $G_2$. 
Following  Artin-Schreier-Witt theory as already used in \cite{MR} (proof of Thm. 5.1, point 2), we introduce the $W_2(\F_p)$-module $$ A:=\frac{ \wp (W_2(L)) \cap W_2(k[X])}{ \wp (W_2(k[X]))}$$ where $W_2(L)$ means the ring of Witt vectors of length $2$ with coordinates in $L$ and $\wp=F-id$.
One can prove that $A$ is isomorphic to the dual of $G_2$ with respect to the Artin-Schreier-Witt pairing (cf. \cite{Bo} Chap. IX, ex. 19). 
Moreover, as a $\Z$-module, $A$ is generated by the classes mod $\wp(k[X])$ of $(f_0(X),g_0(X))$ and $\{(0,f_i(X))\}_{1\leq i \leq r}$  in $W_2(k[X])$. In other words,
$L=k(X,W_i,V_0)_{0 \leq i \leq r}$  is parametrized by the following system of Artin-Schreier-Witt equations:
$$\wp([W_0,V_0])=[f_0(X),g_0(X)] \in W_2(k[X])$$
and
$$\forall \, i \in \{1,\cdots, r\}, \quad \wp(W_i)=f_i(X) \in k[X]$$
An exercise left to the reader shows that one can choose $g_0(X)$ and each $f_i(X)$, with $0 \leq i \leq r$, reduced mod $\wp(k[X])$.

\item \emph{We prove that $f_0 \in \Sigma_2$.} \\
As a $\Z$-module, $p\, A$ is generated by the class of $(0,f_0(X))$ in $A$. By the Artin-Schreier-Witt pairing, $p\,A$ corresponds to the kernel $G_2[p]$ of the map:
 $$ 
\left\{
\begin{array}{ll}
 G_2 \rightarrow G_2\\
g \rightarrow g^p
\end{array}
\right.
$$ 
Thus, $G_2[p] \subsetneq G_2$ is a normal subgroup of $G$. Then, it follows from \cite{MR} (Lemma 2.4) that the pair
$(C/G_2[p], G/G_2[p])$ is a big action parametrized by $W^p-W=f_0(X)$ and with second ramification group $G_2/G_2[p] \simeq \Z/p\Z$. Then, $f_0(X)=X\, S(X)+c\, X \in k[X]$ (cf. \cite{MR} Prop. 2.5), where $S$ is an additive polynomial of $k\{F\}$ with degree $s\geq 1$ in $F$.

\item \emph{The embedding problem.} \\
For any $y \in V$, the classes mod $\wp(k[X])$ of 
$ (f_0(X+y),g_0(X+y))$ and $\{(0,f_i(X+y))\}_{1 \leq i \leq r}$ induces a new generating system of $A$. As in \cite{MR} (proof of Thm 5.1, point 3), this is expressed by the following equation: 
\begin{equation} \label{e9}
\forall \, y \in V, \quad (f_0(X+y),g_0(X+y))= (f_0(X), g_0(X)+ \sum_{i=0}^r \ell_i(y) \,f_i(X)) \, \mod \wp(W_2(k[X]))
\end{equation}
where, for all $i$ in $\{0, \cdots,r\}$, $\ell_i$ is a linear form from $V$ to $\F_p$.
On the second coordinate, \eqref{e9} reads:
\begin{equation} \label{deltag0}
\forall \, y \in V, \quad \Delta_y (g_0):= g_0(X+y)-g_0(X)=\sum_{i=0}^r\,  \ell_i(y) \, f_i(X)+ c \qquad \mod \, \wp (k[X]) 
\end{equation}
where
\begin{equation} \label{eval}
 c =\sum_{i=1}^{p-1}\,  \frac{(-1)^i}{i} \, y^{p-i} \, X^{i+p^{s+1}} +\mbox {lower  degree terms in  X }
\end{equation}
For more details on calculation, we refer to \cite{MR} (proof of Thm 5.1, point 3 and Lemma 5.2).

\item \emph{We prove that $f_i$ lies in $\Sigma_2$, for all $i$ in $\{0,\cdots,r\}$, if and only if $Fratt(G')=[G',G]$.}\\
Put $F:=Fratt(G')$. We deduce from Corollary \ref{Frattbis} that $(C/F,G/F)$ is a big action  whose second ramification group $(G/F)'=(G/F)_2=G_2/F$ is $p$-elementary abelian.
The function field of the curve $C/F$ is now parametrized by the Artin-Schreier equations:
$$\forall \, i \in \{0,\cdots, r\}, \quad \wp(W_i)=f_i(X) \in k[X]$$ 
As $F \subset [G',G]$ (cf. Lemma \ref{Frattbis}), 
$$F=[G',G]=[G_2:G] \Leftrightarrow \{e\} =[G_2/F,G/F]=[(G/F)',G/F]\Leftrightarrow (G/F)'\subset Z(G/F)$$
By \cite{MR2} (Prop.2.13), this occurs if and only if for all $ i$ in $\{0,\cdots, r\}$, $f_i(X)=X\, S_i(X)+c_i\, X \in \Sigma_2$.

\item \textit{We prove that $g_0$ does not belong to $\Sigma_p$.}\\
We first notice that the right-hand side of \eqref{deltag0} does not belong to $\Sigma_{p-1}$: indeed, the monomial $X^{p-1+p^{s+1}} \in \Sigma_p-\Sigma_{p-1}$ occurs once in $c$ but not in $\sum_{i=0}^r\,  \ell_i(y) \, f_i(X)$ which lies in $\Sigma_2 \subset \Sigma_{p-1}$, for $p \geq 3$. Now, assume that $g_0 \in \Sigma_p$. Then, by \cite{MR2} (Lemma 3.9), the left-hand side of \eqref{deltag0}, namely $\Delta_y (g_0)$, lies in $\Sigma_{p-1}$, hence a contradiction.
Therefore, one can define an integer $a$ such that $X^{a}$ is the monomial of $g_0(X)$ with highest degree among those that do not belong to $\Sigma_p$. Note that since $g_0$ is reduced mod $\wp(k[X])$, $a \neq 0$ mod $p$.

\item \textit{We prove that $a-1 \geq p-1+p^{s+1}$.}\\
We have already seen that the monomial $X^{p-1+p^{s+1}}$ occurs in the right hand side of \eqref{deltag0}.
In the left-hand side of \eqref{deltag0}, $X^{p-1+p^{s+1}}$ is produced by monomials $X^b$ of $g_0$ with $b>p-1+p^{s+1}$.
If $b >a$, $X^b \in \Sigma_p$, so $\Delta_y(X^b) \in \Sigma_{p-1}$, which is not the case of $X^{p-1+p^{s+1}}$.
It follows that $X^{p-1+p^{s+1}}$ comes from monomials $X^b$ with $a\geq b > p-1+p^{s+1}$. Hence the expected inequality.

\item \textit{We prove that $p$ divides $a-1$.}\\
Assume that $p$ does not divide $a-1$. In this case, the monomial $X^{a-1}$ is reduced mod $\wp(k[X])$ and \eqref{deltag0} reads as follows:
$$\forall \, y \in V, \quad c_a(g_0) \, a\, y \, X^{a-1}+ S_{p-1}(X)+ R_{a-2}(X)= c +\sum_{i=0}^r \, \ell_i(y) \, f_i(X) \,  \mod \, \wp\, (k[X]) $$ where $c_a(g_0)\neq 0$ denotes the coefficient of $X^a$ in $g_0$, $S_{p-1}(X)$ is a polynomial in $\Sigma_{p-1}$ produced by  monomials $X^b$ of $g_0$ with $b > a$ and $R_{a-2}(X)$ is a polynomial of $k[X]$ with degree lower than $a-2$  produced by  monomials $X^b$ of $g_0$ with $b \leq a$.
We first notice that $X^{a-1}$ does not occur in $S_{p-1}(X)$. Otherwise, $X^{a-1} \in \Sigma_{p-1}$ and $X^a =X^{a-1} \, X \in \Sigma_p$, hence a contradiction. Likewise, $X^{a-1}$ does not occur in $\sum_{i=0}^r \, \ell_i(y) \, f_i(X) \in \Sigma_2$. Otherwise, $X^{a}=X^{a-1} \, X \in \Sigma_3 \subset \Sigma_p$, as $p \geq 3$. It follows that $X^{a-1}$ occurs in $c$, which requires $a-1 \leq$ deg $b=p-1+p^{s+1}$. Then, the previous point implies $a-1 = p-1+p^{s+1}$, which contradicts $a \neq 0$ mod $p$.\\
Thus, $p$ divides $a-1$. So, we can write $a=1+\lambda\, p^t$, with $t>0$, $\lambda$ prime to $p$ and $\lambda \geq 2$ because of the definition of $a$. We also define $j_0:=a-p^t=1+(\lambda-1)\, p^t$. 

\item \textit{We search for the coefficient of the monomial $X^{j_0}$ in the left-hand side of \eqref{deltag0}. }\\
Since $p$ does not divide $j_0$, the monomial $X^{j_0}$ is reduced mod $\wp(k[X])$.
 In the left-hand side of \eqref{deltag0}, namely $\Delta_y(g_0)$ mod $\wp(k[X])$, the monomial 
$X^{j_0}$ comes from monomials of $g_0(X)$ of the form: $X^b$, with $b \geq j_0+1$. However, as seen above, the monomials $X^b$ with $b>a$ produce in $\Delta_y(g_0)$ elements that belong to $\Sigma_ {p-1}$, whereas $X^{j_0} \not \in \Sigma_{p-1}$. Otherwise, $X^a=X^{j_0} \, X^{p^t} \in \Sigma_p$, which contradicts the definition of $a$. So we only have to consider the monomials $X^b$ of $g_0(X)$ with $b \in \{j_0+1,\cdots, a\}$. Then, the same arguments as those used in \cite{MR} (proof of Thm. 5.1, point 6) allow to conclude that 
the coefficient of $X^{j_0}$ in the left-hand side of \eqref{deltag0} is $T(y)$ where $T(Y)$ denotes a polynomial of $k[X]$ with degree $p^t$. 

\item \emph{We identify with the coefficient of $X^{j_0}$ in the right-hand side of \eqref{deltag0} and gather a contradiction.}
 As mentionned above, the monomial $X^{j_0}$ does not occur in $\sum_{i=0}^r\,  \ell_i(y) \, f_i(X) \in \Sigma_2 \subset \Sigma_{p-1}$, for $p \geq 3$. 
Assume that the monomial $X^{j_0}$ appears in $c$, which implies that $j_0 \leq p-1+p^{s+1}$. Using the same arguments as in \cite{MR} (proof of Thm. 5.1, point 7), we gather that $j_0=1+(\lambda-1)\, p^t =1+p^{s+1}$. Then, $X^{j_0}$ lies in
 $\Sigma_2$, which leads to the same contradiction as above.
Therefore, the monomial $X^{j_0}$ does not occur in the right-hand side of \eqref{deltag0}. Then, $T(y)=0$ for all $y$ in $V$, which means that $|V| \leq p^t$. 
Call $C_0$ the curve whose function field is parametrized by $\wp([W_0,V_0])=[f_0(X),g_0(X)]$. 
The same calculation as in \cite{MR} (proof of Theorem 5.1, point 7) shows that $g_{C_0} \geq p^{t+1}\, (p-1)$.
Furthermore, $g \geq p^r \,g_{C_0}$
 (see e.g. \cite{LM} Prop. 8.5, formula (8)).
As $|G|=|G_2||V| \leq p^{2+r+t}$, it follows that $\frac{|G|}{g}=\frac{p}{p-1} <\frac{2\,p}{p-1}$, hence a contradiction.
 $\square$
\end{enumerate}

\section{Classification of big actions under condition $(*)$.}

We now pursue the classification of big actions initiated by Lehr and Matignon who characterize big actions $(C,G)$ satisfying $\frac{|G|}{g^2} \geq \frac{4}{(p-1)^2}$ (cf. \cite{LM}). In this section, we exhibit a parametrization for big actions $(C,G)$ satisfying condition $(*)$, namely: 
$$\frac{|G|}{g^2} \geq \frac{4}{(p^2-1)^2} \qquad \qquad  (*)$$ 
As proved in \cite{MR} (Prop. 4.1  and Prop. 4.2), this condition requires $G'(=G_2)$ to be an elementary abelian $p$-group  with order dividing $p^3$. Since $G_2$ cannot be trivial (cf. \cite{MR} Prop. 2.2), this leaves three possibilities.
This motivates the following
\begin{defi} \label{def11}
Let $(C,G)$ be abig action. Let $i \geq 1$ be an integer.
We say that
\begin{enumerate}
\item $(C,G)$ satisfies $\mathcal{G}_*$ if $(C,G)$ satisfies condition $(*)$
\item $(C,G)$ satisfies $\mathcal{G}_*^{p^i}$ if $(C,G)$ satisfies $\mathcal{G}_*$ with $G'\simeq (\Z/p\Z)^i$.
\end{enumerate}
\end{defi}

\subsection{Preliminaries: big actions with a $p$-elementary abelian $G'(=G_2)$.}
To start with, we fix the notations and recall some necessary results on big actions with a $p$-elementary abelian $G_2$ drawn from \cite{MR2}.

\medskip

\begin{recall}
Let $(C,G)$ be a big action such that $G'(=G_2) \simeq (\Z/p\Z)^n$, $n \geq 1$. Write the exact sequence:
$$0 \longrightarrow G_2 \simeq (\Z/p\Z)^n
\longrightarrow G\stackrel{\pi}{\longrightarrow} V \simeq (\Z/p\Z)^v \longrightarrow 0 
$$
\begin{enumerate}
\item  We denote by $L$ be the function field of the curve $C$ and by $k(X):=L^{G_2}$ the subfield of $L$ fixed by $G_2$.
Then, the extension $L/k(X)$ can be parametrized by $n$ Artin-Schreier equations: $W_i^p-W_i=f_i(X) \in k[X]$ with $1\leq i \leq n$. Following \cite{MR2} (Def. 2.3), one can choose an "adapted basis" $\{f_1(X), \cdots, f_n(X)\}$ with some specific properties:
\begin{enumerate}
\item  For all $i \in \{1,\cdots,n\}$, each function $f_i$ is assumed to be reduced mod $\wp(k[X])$ 
\item For all $i \in \{1,\cdots,n\}$, put $m_i:=deg \,f_i$. Then, $m_1 \leq m_2 \leq \cdots \leq m_n$.
\item  $\forall \;  (\lambda_1,\cdots \lambda_n) \in \F_p^n$ not all zeros, 
 $$ deg \; ( \sum_{i=1}^n \, \lambda_i  \, f_i(X)) = \max_{i\in \{1,\cdots,n\}} \{deg \;  \lambda_i \,f_i(X)  \}.$$
 \end{enumerate}
In this case, the genus of the curve $C$ is given by the following formula (cf. \cite{MR2} Cor. 2.7):
\begin{equation} \label{genus}
g=\frac{p-1}{2} \sum_{i=1}^n \, p^{i-1} \, (m_i-1)
\end{equation}
\item Now, consider the $\F_p$-subvector space of $k[X]$ generated by the classes of $\{f_1(X),\cdots, f_n(X)\}$ mod $\wp(k[X])$:
$$A:= \frac{ \wp (L) \cap k[X]}{ \wp (k[X])} $$
Recall that $A$ is isomorphic to the dual of $G_2$ with respect to the Artin-Schreier pairing (cf. \cite{MR2} section 2.1). As seen in \cite{MR2} (section 2.2), $V$ acts on $G_2$ via conjugation. This induces a representation $\phi$: $V \rightarrow Aut(G_2)$. The representation $\rho: V \rightarrow Aut(A)$, which is dual with respect to the Artin-Schreier pairing, expresses the action of $V$ on $A$ by translation. More precisely, for all $y$ in $V$, the automorphism $\rho(y)$ is defined as follows:
 $$  \rho(y): 
\left\{
\begin{array}{lc}
A \rightarrow A \\
\overline{f(X)} \rightarrow \overline{f(X+y)}
\end{array}
\right.
$$ where $\overline{f(X)}$ means the class in $A$ of $f(X) \in k[X]$ 
For all $y$ in $V$, the matrix of the automorphism $\rho(y)$ in the adapted basis fixed for $A$ is an upper triangular matrix of $Gl_n(\F_p)$ with identity on the diagonal, namely
 $$ L(y):=
\begin{pmatrix}
 1 &\ell_{1,2}(y)&  \ell_{1,3}(y) &\cdots& \ell_{1,n}(y) \\
 0 & 1 &\ell_{2,3}(y) & \cdots&\ell_{2,n}(y) \\
 0 & 0 &\cdots&\cdots&\ell_{i,n}(y) \\
 0 & 0 & 0 & 1& \ell_{n-1,n}(y)\\
 0 & 0 & 0 & 0 & 1
\end{pmatrix} \in Gl_n(\F_p)
$$
where, for all $i$ in $\{1,\cdots, n-1\}$, $\ell_{i,i+1}$ is a nonzero linear form from $V$ to $\F_p$ (see \cite{MR2} section 2.4). In other words, 
$$ \forall \, y \in V, \; f_1(X+y)-f_1(X)=0 \quad \mod \; \wp(k[X])$$
\begin{equation} \label{toto}
\forall \, i \in \{2,\cdots,n\}, \; \forall \, y \in V, \; f_i(X+y)-f_i(X)=\sum_{j=1}^{i-1}\,
\ell_{j,i}(y) \, f_j(X) \quad \mod \; \wp(k[X]) 
\end{equation}
For all map $\ell$, we write $\ell=0$ if $\ell$ is identically zero and $\ell \neq 0$ otherwise.
\item The case of a trivial representation can be charactrized as follows (see \cite{MR2} Prop. 2.13).
Indeed, the following assertions are equivalent:
\begin{enumerate}
\item The representation $\rho$ is trivial, i.e. $$\forall \, i \in \{1,\cdots,n\}, \quad \forall \, y \in V, \quad f_i(X+y)-f_i(X)=0 \mod \; \wp(k[X])$$
\item The commutator subgroup of $G'$ and $G$ is trivial, i.e.  $G' \subset Z(G)$.
\item For all $i$ in $\{1,\cdots,n\}$, $f_i(X)=X \, S_i(X)+ c_i\,X \in k[X]$  where each $S_i \in k\{F\}$ is an additive polynomial with degree $s_i \geq 1$ in $F$. 
 So, write $ S_i(F) = \sum _{j=0}^{s_i} \,a_{i,j} \,F^j$ with $a_{i,s_i} \neq 0$. Then, one defines an additive polynomial related to $f_i$, called the "palindromic polynomial" of $f_i$:
 $$Ad_{f_i}:= \frac{1}{a_{i,s_i}^{p^{s_i}}} \; F^{s_i}  (\sum _{j=0}^ {s_i} \,a_{i,j} \,F^j + F^{-j} \,a_{i,j})$$
In this case, $$V \subset \bigcap_{i=1}^n Z (Ad_{f_i}) $$
\end{enumerate}
Since, under condition $(*)$, $G'$ is $p$-elementary abelian, we deduce from point $(b)$ that the case of a trivial representation corresponds to the case $\{e\}=Fratt(G')=[G',G]$.
\item To conclude, we recall that for all $t \geq 1$, $\Sigma_t$ means the $k$-subvector space of $k[X]$ generated by $1$ and the products of at most $t$ additive polynomials of $k[X]$ (cf. \cite{MR2} Def. 3.1). As proved in \cite{MR2} (Thm. 3.13), for all $i$ in $\{1,\cdots, n\}$, $f_i$ lies in $\Sigma_{i+1}$.
\end{enumerate}
\end{recall}

\subsection{First case: big actions satisfying $\mathcal{G}_*^p$.}

\begin{proposition} \label{caspcyclique}
\begin{enumerate}
We keep the notations of section 5.1.
\item $(C,G)$ is a big action with $G_2\simeq \Z/p\Z$ if and only if $C$ is birational to a curve $C_f$ parametrized by $W^p-W=f(X)=X\, S(X)\in k[X]$, where $S$ is a (monic) additive polynomial with degree $s\geq 1$ in $F$. 
\item In what follows, we assume that $C$ is birational to a curve $C_f$ as described in the first point. 
\begin{enumerate}
\item If $s\geq 2$, $A_{\infty,1}$ is the unique $p$-Sylow subgroup of $A$, where $\infty$ denotes the point of $C$ corresponding to $X=\infty$.
\item If $s=1$, there exists $r:=p^3+1$ points of $C$: $P_0:=\infty, P_1, \cdots,P_r$ such that 
$(A_{P_i,1})_{0 \leq i \leq r}$ are the $p$-Sylow subgroups of $A$. In this case, for all $i$ in $\{1,\cdots,r\}$, there exists $\sigma_i \in A$ such that $\sigma_i(P_i)=\infty$.
\end{enumerate}
In both cases, $A_{\infty,1}$ is an extraspecial group (see \cite{Su2} Def. 4.14) with exponent $p$ (resp. $p^2$) if $p>2$ (resp. $p=2$) and order $p^{2s+1}$.
More precisely, $A_{\infty,1}$ is a central extension of its center $Z(A_{\infty,1})=(A_{\infty,1})'$ by the elementary abelian $p$-group $Z(Ad_f)$, i.e.
$$ 0 \longrightarrow Z(A_{\infty,1})=(A_{\infty,1})'\simeq \Z/p \Z \longrightarrow A_{\infty,1}  \stackrel{\pi}{\longrightarrow}  Z(Ad_f) \simeq (\Z / p \Z) ^{2s} \longrightarrow 0$$ 
Furthermore, $(C,A_{\infty,1})$, and so each $(C,A_{P_i,1})$, with $1 \leq i \leq r$, are big actions satisfying $\mathcal{G}_*^p$.
\item Let $V$ be a subvector space of $Z(Ad_f)$ with dimension $v$ over $\F_p$. Then, $(C,\pi^{-1}(V))$ is also a big action satisfying $\mathcal{G}_*^p$ if and only if

$$
\begin{array}{ll} \label{conditions}
 \mbox{if} \, \, p \neq 2, \, & 2s \geq v \geq \max \{s+1,2s-3\}\\
 \mbox{if} \, \, p = 2, \, & 2s \geq v \geq \max \{s+1,2s-4\}
\end{array}
$$
\end{enumerate}
\end{proposition}

\noindent We collect the different possibilities in the table below:\\

\medskip

\noindent \begin{tabular}{|c|c|c|c|c|}
\hline
case & $v$ & $s$ & $V$ &  $G$  \\
\hline
\hline
1-& $2s$ & $s\geq 1^{\dagger}$ & $Z(Ad_f)^{\dagger}$ & $A_{\infty,1}^{\dagger}$ \\
\hline
2 & $2s-1$ &  $s\geq 2$  & index p subgroup of $Z(Ad_f)$&  index p subgroup of $A_{\infty,1}$ \\
\hline
3 & $2s-2$ & $s\geq 3$  &index $p^2$ subgroup of $Z(Ad_f)$& index $p^2$ subgroup of $A_{\infty,1}$ \\
\hline
4 & $2s-3$ & $s\geq 4$  &index $p^3$ subgroup of $Z(Ad_f)$&  index $p^3$ subgroup of $A_{\infty,1}$ \\
\hline
5 (p=2)& $2s-4$ & $s\geq 5$  &index $p^4$ subgroup of $Z(Ad_f)$&  index $p^4$ subgroup of $A_{\infty,1}$ \\
\hline
\end{tabular}

\bigskip

\noindent \begin{tabular}{|c|c|c|}
\hline
case & $|G|/ g$ & $|G| /g^2$ \\
\hline
\hline
1  & $\frac{2\,p}{p-1}\,  p^{s}$ &$\frac{4}{(p^2-1)^2} \, (p+1)^2 \, p $\\
\hline
2  &$\frac{2\,p}{p-1} \,  p^{s-1}$ & $\frac{4}{(p^2-1)^2} \, (p+1)^2 $\\
\hline
3  &$\frac{2\,p}{p-1}\, p^{s-2}$ & $\frac{4}{(p^2-1)^2} \, \frac{(p+1)^2}{p}$\\
\hline
4  &$\frac{2\,p}{p-1}\,  p^{s-3}$ & $\frac{4}{(p^2-1)^2} \, \frac{(p+1)^2}{p^2}$\\
\hline
5 (p=2)& $\frac{2\,p}{p-1}\,  p^{s-4}$ & $\frac{4}{(p^2-1)^2} \, \frac{(p+1)^2}{p^3}$\\
\hline
\end{tabular}

\bigskip

\noindent $\dagger$ \textit{Note:} In the case $s=1$, this result is true up to conjugation by $\sigma_i$ as defined in Proposition \ref{caspcyclique}.

\medskip

\noindent \textbf{Proof:}
\begin{enumerate}
\item See \cite{LM} (Thm. 1.1 I)
\item See Remark \ref{italien},  \cite{LM} (Thm. 3.1) and \cite{MR} (Prop. 2.5).
\item  This essentially derives from Proposition \ref{prop11} which implies $(p+1)^2 \geq p^{2s-v-1}$. If $2s-v-1 \geq 3$, it implies $p^2+2\,p+1 \geq p^3$, which is impossible for $p>2$. Accordingly, if $p>2$, we obtain $2s-v-1 \leq 2$, which means $v \geq 2s-3$. If $p=2$, $ (p+1)^2 \geq p^{2s-v-1}$ is satisfied if and only if $2s-v-1 \leq 3$, i.e. $v \geq 2s-4$. The claim follows. $\square$
\end{enumerate}
\medskip

\begin{remarque}
Note that, for $p>2$, the solutions can be parametrized by $s$ algebraically independent variables over $\F_p$, namely the $s$ coefficients of $S$ assumed monic after an homothety on the variable $X$. Note that $s \sim log \, g$.
\end{remarque}

\subsection{Second case: big actions satisfying $\mathcal{G}_*^{p^2}$.}

\subsubsection{Case: $[G',G]=Fratt(G')=\{e\}$.}

\begin{proposition}  \label{casp2trivial}
Let $(C,G)$ be a big action satisfying $\mathcal{G}_*^{p^2}$. Assume that $[G',G]=\{e\}$ and keep the notations of section 5.1.
\begin{enumerate}
\item The pair $(C, A_{\infty,1})$ is a big action satisfying $\mathcal{G}_*^{p^2}$ 
. Moreover, $A_{\infty,1}$ is a special group (see \cite{Su2} Def. 4.14) with exponent $p$ (resp. $p^2$)
   (for $p>2$ (resp. $p=2$) and order $p^{2+2\, s_1}$.
More precisely, $A_{\infty,1}$ is a central extension of its center $Z(A_{\infty,1})=(A_{\infty,1})'$ by the elementary abelian $p$-group $Z(Ad_{f_1})$, i.e.
$$ 0 \longrightarrow Z(A_{\infty,1})=(A_{\infty,1})'\simeq (\Z/p \Z)^2 \longrightarrow A_{\infty,1}  \stackrel{\pi}{\longrightarrow}  Z(Ad_{f_1}) \simeq (\Z / p \Z) ^{2s_1} \longrightarrow 0$$ 
\item Furthermore, $s_2=s_1$ or $s_2=s_1+1$.
\begin{enumerate}
\item If $s_2=s_1$, 
$G=\pi^{-1}(V)$, where $V$ is a subvector space of $Z(Ad_{f_1})$ with dimension $v$ over $\F_p$ such that 
$2\, s_1-2 \leq v \leq 2\,s_1$. Then, $A_{\infty,1}$ is a $p$-Sylow subgroup of $A$. It is normal except if $C$ is birationnal to the Hermitian curve: $W^q-W=X^{1+q}$ with $q=p^2$.
\item If $s_2=s_1+1$, $V=Z(Ad_{f_1})$ and $G=A_{\infty,1}$ is the unique $p$-Sylow subgroup of $A$.
\end{enumerate}
\end{enumerate}
\end{proposition}

The different possibilities are listed in the table below:\\
\medskip

\noindent \begin{tabular}{|c|c|c|c|c|c|}
\hline
case & $s_1$ & $s_2$ & $v$ & $V$ & $G$   \\
\hline
\hline
(a)-1 & $s\geq 2$ & $s$ & $2 s$ & $Z(Ad_{f_1})=Z(Ad_{f_2})$ &$A_{\infty,1}$ \\
\hline
(a)-2 & $s\geq 2$ & $s$ & $2 s -1 $ & index $p$ subgroup of $Z(Ad_{f_1})$ & index $p$ subgroup of $A_{\infty,1}$\\
\hline
(a)-3&  $s\geq 3$ & $s$ & $2 s -2 $ & index $p^2$ subgroup of $Z(Ad_{f_1})$ & index $p^2$ subgroup of $A_{\infty,1}$\\
\hline
\hline
(b) & $s\geq 3$ & $s+1$ & $2 s$ & $Z(Ad_{f_1})$ &$A_{\infty,1}$\\
\hline
\end{tabular}

\bigskip

\noindent \begin{tabular}{|c|c|c|}
\hline
case & $ |G| / g$ & $ |G| /g^2$  \\
\hline
\hline
(a)-1&  $\frac{2\,p}{p-1} \, \frac{p^{1+s}}{1+p}$ & $\frac{4}{(p^2-1)^2} \,  p^2$\\
\hline
(a)-2& $\frac{2\,p}{p-1} \, \frac{p^{s}}{1+p}$ & $\frac{4}{(p^2-1)^2} \, p$ \\
\hline
(a)-3& $\frac{2\,p}{p-1} \, \frac{p^{s-1}}{1+p}$ & $\frac{4}{(p^2-1)^2} $ \\
\hline
\hline
(b)& $\frac{2\,p}{p-1} \, \frac{p^{1+s}}{1+p}$ & $\frac{4}{(p^2-1)^2} \, \frac{p^2(p+1)^2}{(1+p^2)^2}$ \\
\hline
\end{tabular}

\bigskip

\noindent \textbf{Proof:} 
\begin{enumerate}
\item Use Proposition \ref{prop11} to prove that the pair $(C,A_{\infty,1})$ is a big action satisfying $\mathcal{G}_*^{p^2}$
with the following exact sequence:
 $$ 0 \longrightarrow A_{\infty,2} \longrightarrow A_{\infty,1}  \stackrel{\pi}{\longrightarrow} Z(Ad_{f_1}) \simeq
   (\Z/\,p\, \Z)^{2\,s_1} \longrightarrow 0$$
The proof to show that $A_{\infty,1}$ is a special group, i.e. satisfies $Z(A_{\infty,1})=(A_{\infty,1})'=Fratt(A_{\infty,1})\simeq (\Z/p\Z)^2$, is the same that the one exposed in \cite{MR2} (Prop. 4.3.3). Nevertheless, one has to choose $\mathcal{H}$ an index $p$-subgroup of $G_2$ such that $C/\mathcal{H}$ is the curve parametrized by $W_1^p-W_1=f_1(X)$.
\item Assume that $s_2-s_1 \geq 2$. Then, $|G|=p^{2+v} \leq p^{2+2\,s_1}$ and $g=\frac{p-1}{2} \, p^{s_1} \, (1+p^{1+s_2-s_1}) \geq \frac{p-1}{2} \, p^{s_1} \, (1+p^{3})$.
So, $\frac{|G|}{g^2} \leq \frac{4}{(p^2-1)^2} \, \frac{(1+p)^2\, p^2}{(1+p^3)^2} <\frac{4}{(p^2-1)^2}$, which contradicts condition $(*)$. So, $0 \leq s_2-s_1 \leq 1$. In each case, the description of $A_{\infty,1}$ and $G$ derive from Proposition \ref{prop11} combined with Remark \ref{italien}. $\square$
\end{enumerate}

\medskip

To go further in the description of the functions $f_i's$ in each case, we introduce two additive polynomials
$\mathcal{V}$ and $T$ defined as follows:
$$\forall \, i \in \{1,2\}, \quad \mathcal{V}:= \prod_{y \in V}\,  (X-y)\quad  \mbox{divides} \quad T:= gcd \{ Ad_{f_1}, Ad_{f_2} \} \quad \mbox{divides}
\quad  Ad_{f_i}$$
\noindent In what follows, we work in the Ore ring $k\{F\}$ and write the additive polynomials as polynomials in $F$.

\medskip

\noindent \begin{tabular}{|c|c|c|c|c|c|}
\hline
case & $deg_F \, \mathcal{V}$ & $deg_F \, T$ & $deg_F \,(Ad_{f_1})$ &$deg_F \,(Ad_{f_2}) $ &$T$ \\
\hline
\hline
(a)-1 & $2s$ & $2s$ & $2s$ & $2s$ & $\mathcal{V} = T=Ad_{f_1}= Ad_{f_2}$ \\
\hline
\hline
(a)-2-i& $2s-1$ & $2s$ & $2s$ & $2s$ & $\mathcal{V} \quad \mbox{divides} \quad T=Ad_{f_1}= Ad_{f_2}$ \\
\hline
(a)-2-ii & $2s-1$ &$2s-1$ & $2s$ & $2s$ & $\mathcal{V}=T \quad \mbox{divides} \quad Ad_{f_1} $ \\
\hline
\hline
(a)-3-i & $2s-2$ & $2s$ & $2s$ & $2s$ & $\mathcal{V} \quad \mbox{divides} \quad T=Ad_{f_1}= Ad_{f_2}$ \\
\hline
(a)-3-ii & $2s-2$ & $2s-1$  & $2s$ & $2s$& $\mathcal{V} \quad \mbox{divides} \quad T \quad \mbox{divides} \quad Ad_{f_1}$ \\
 \hline
(a)-3-iii& $2s-2$  & $2s-2$ & $2s$ & $2s$&  $\mathcal{V}=T \quad \mbox{divides} \quad Ad_{f_1} $\\
\hline
\hline
(b) & $2s$ & $2s$ & $2s$ & $2s+2$ & $\mathcal{V}=T=Ad_{f_1} \quad \mbox{divides} \quad Ad_{f_2} $ \\
\hline 
\end{tabular}

\bigskip

\noindent The three cases where $Ad_{f_1}=Ad_{f_2}$ can be parametrized in the same way as in \cite{MR2} (Prop. 4.2).

\bigskip

\noindent \begin{tabular}{ | c | c | c | }
\hline
case &   $S_1$ or $Ad_{f_1} $ & $S_2$ or $Ad_{f_2}$\\
\hline
\hline
(a)-1 &  $S_1=\sum_{j=0}^{s/d} \, \alpha_{jd} \, F^{jd}, \quad \alpha_s=1 $  &  $S_2=\gamma \, S_1, \quad  \gamma \in \F_{p^d} - \F_p, \quad d \geq 2$ \\
\hline
\hline
(a)-2-i & $S_1 =\sum_{j=0}^{s/d} \, \alpha_{jd} \, F^{jd},\quad \alpha_s=1 $  &  $S_2=\gamma \,S_1, \quad  \gamma \in \F_{p^d} - \F_p,\quad d \geq 2 $\\
\hline
(a)-2-ii &  $Ad_{f_1}= (\alpha_1 \,F +\beta_1 \,I) \, T, \quad \alpha_1 \neq 0$&    $Ad_{f_2}=(\alpha_2 \,F +\beta_2 \,I)\, T, \quad \alpha_2 \neq 0 $\\
\hline
\hline
(a)-3-i  & $S_1 =\sum_{j=0}^{s/d}\,  \alpha_{jd} \,F^{jd}, \quad \alpha_s=1 $  &  $S_2=\gamma \, S_1, \quad  \gamma \in \F_{p^d} - \F_p, \quad d \geq 2$\\
\hline
(a)-3-ii &  $Ad_{f_1}= (\alpha_1\, F +\beta_1\, I) \, T, \quad \alpha_1 \neq 0$&     $Ad_{f_2}=(\alpha_2\, F +\beta_2\, I)\, T, \quad \alpha_2 \neq 0 $\\
 \hline
(a)-3-iii   & $Ad_{f_1}= (\alpha_1\, F^2 +\beta_1 \,F +\delta_1\, I ) \, T,\quad \alpha_1 \neq 0$&     $Ad_{f_2}=(\alpha_2 \,F^2 +\beta_2\, F+\delta_{2}\, I )\, T , \quad \alpha_2 \neq 0$\\
\hline
\hline
(b) &$ Ad_{f_1}=\prod_{v \in V} (X-v)$  & $Ad_{f_2} = (\alpha_2 \,F^2 + \beta_2 \,F +\delta_2\, I) \, Ad_{f_1}, \quad \alpha_2 \neq 0$ \\
\hline 
\end{tabular}

\bigskip

\noindent We display the parametrization of the functions $f_i$'s in the case (a)-2-ii for the smallest values of $s$, namely $s=2$ and $s=3$. 

\bigskip

\noindent \textbf{Cas (a)-2-ii with $s=2$ for $p>2$.}

\bigskip

\noindent \begin{tabular}{ | c | c | }
\hline
$f_1$ & $f_1(X)=X^{1+p^2} + a_{1+p} \, X^{1+p} + \frac{1}{2} \, a_2\,  X^2$\\
\hline
$a_{1+p}$ & $a_{1+p} \in k$ \\
\hline
$a_{2}$ & $a_{2} \in k$ \\
\hline
\hline
$f_2$ & $f_2(X)= b_{1+p^2}^{p^2}\,  X^{1+p^2} + b_{1+p} \,  X^{1+p}+ b_2\, X^2 +b_{1}\, X$\\
\hline
$b_{1+p^2}$ & $ b_{1+p^2} \in Z (w^{p^2} \, X^{p^3} +w^p _,(-a_{2}^p +a_{1+p}^p\, w^{p^2} -w^{p^2+p^3}) \, X^{p^2}  +(a_{1+p}-w^{p^2} ) \, X^p -w^{-1} \, X)$ \\
& \mbox{with} $\quad b_{1+p^2} \not \in \F_{p^{2}}$.\\
\hline
$w$ & $w \in Z (X^{1+p+p^2+p^3} -a_{1+p}^p \, X^{1+p+p^2} +a_{2}^p\,  X^{1+p} -a_{1+p} \, X+1)$\\
\hline
$b_{1+p}$ & $b_{1+p}=w^{p^2} (b_{1+p^2}^{p^2}-b_{1+p^2})^p +b_{1+p^2}^p \, a_{1+p} $\\
\hline 
$b_2$ & $2\, b_2=w^p \,(b_{1+p^2}^{p^2}-b_{1+p^2})\, (a_{1+p}-w^{p^2})+b_{1+p^2}\,_{2}$ \\
\hline 
$b_1$ & $b_1\in k$\\
\hline
\end{tabular}

\bigskip

\noindent \textbf{Case (a)-2-ii with $s=3$ for $p>2$.}

\bigskip

\noindent \begin{tabular}{ | c | c | }
\hline
$f_1$ & $f_1(X)=X^{1+p^3} + a_{1+p^2} \, X^{1+p^2} + a_{1+p} \,  X^{1+p} + \frac{1}{2} \, a_2 \, X^2$\\
\hline
$a_{1+p^2}$ & $a_{1+p^2} \in k$ \\
\hline
$a_{1+p}$ & $a_{1+p} \in k$ \\
\hline
$a_{2}$ & $a_{2} \in k$ \\
\hline
\hline
$f_2$ & $f_2(X)= b_{1+p^3}^{p^3}\, X^{1+p^3} +b_{1+p^2}\,  X^{1+p^2} +b_{1+p} \, X^{1+p} +b_2\, X^2 +b_{1} \, X$\\
\hline
$w$ & $w \in Z (X^{1+p+p^2+p^3+p^4+p^5} -a_{1+p^2}^{p^2} X^{1+p+p^2+p^3+p^4} $\\
&$ +a_{1+p}^{p^2} X^{1+p+p^2+p^3} -a_{2}^{p^2} X^{1+p+p^2} +a_{1+p}^p X^{1+p}- a_{1+p^2}X+1)$\\
\hline
$b_{1+p^3}$ & $b_{1+p^3} \in Z( P_1)\cap Z(P_2)-\F_ {p^3} $ \\ 
& \mbox{with} $\quad P_1(X)=w^{p^3+1} \, X^{p^5} + (1-w\, a_{1+p^2}) X^{p^3} +(w\, a_{1+p^2}-w^{p^3+1})\, X^{p^2 }-X$ \\
& \mbox{with} $\quad P_2(X)=w^{p^2} (a_{1+p^2}-w^{p^3}) \, X^{p^4}+ w^p(-a_{2}^p+a_{1+p}^p \, w^{p^2} -a_{1+p^2}^p w^{p^2+p^3} +w^{p^2+p^3+p^4}) \, X^{p^3}$ \\
&$ +(a_{1+p}+w^{p^2+p^3} -a_{1+p^2} w^{p2}) \,X^p $\\
&$+(-a_{1+p}+a_{2}^p w^p-a_{1+p}^p w^{p+p^2} +a_{1+p^2}^p w^{p+p^2+p^3} -w^{p+p^2+p^3+p^4}) \, X$\\
\hline
$b_{1+p^2}$ & $b_{1+p^2}=w^{p^3} \, (b_{1+p^3}^{p^3}-b_{1+p^3})^{p^2} +b_{1+p^3}^{p^2}\, a_{1+p^2}$\\
\hline 
$b_{1+p}$ & $b_{1+p}= w^{p^2}\,  (b_{1+p^3}^{p^3} -b_{1+p^3})^p \, (a_{1+p^2}-w^{p^3})+ b_{1+p^3}^p \, a_{1+p}$\\
\hline
$b_2$ & $2\, b_2= w^p\, (b_{1+p^3}^{p^3}-b_{1+p^3})\, (a_{1+p}-a_{1+p^2}\, w^{p^2} +w^{p^2+p^3}) +b_{1+p^3}\,a_{2} $ \\
\hline 
$b_1$ & $b_1 \in k$\\
\hline
\end{tabular}

\bigskip

\noindent The calculation of the case $s=3$ already raises a problem as the parameter $b_{1+p^3}$ has to lie in the set of zeroes of two polynomials.

\bigskip

\noindent For the remaining last two cases (a)-3-iii and (b), we merely display examples of realization so as to prove the effectiveness of these cases. 

\bigskip

\noindent \textbf{An example of realization for the case (a)-3-iii.}

\medskip

\noindent \begin{tabular}{ | c | c | }
\hline
$T$ & $T =F^{2s -2}+I$\\
\hline
$V$ & $V=Z(F^{2s -2}+I)$\\
\hline
$Ad_{f_1}$ & $Ad_{f_1}= (F^2+I)\,  T $ \\
\hline 
$f_1$ & $f_1(X)=X^{1+p^{s}} + X^{1+p^{s -2}}$\\
\hline
\hline
$Ad_{f_2} $ & $ Ad_{f_2}= (F^2+F+I)\, T $ \\
\hline 
$f_2$ & $f_2(X)= X^{1+p^{s}} + X^{1+p^{s -1}}+ X^{1+p^{s-2}}$\\
\hline
\end{tabular}

\bigskip

\noindent \textbf{An example of realization for the case (b).}

\medskip

\noindent \begin{tabular}{ | c | c | }
\hline
$f_1$ & $f_1(X)=X^{1+p^{s}}$\\
\hline
\hline 
$f_2$ & $f_2(X)= \alpha_2 \, X^{1+p^{s +1}} + \beta_2\,  X^{1+p^{s}} + \delta_2 \, X^{1+p^{s -1}}$\\
\hline
$\alpha_2$ & $\alpha_2 \in \F_{p^{2 s}}$ \\
\hline
$\beta_2$ & $\beta_2 \in \F_{p^{s}}$\\
\hline
$\delta_2$ & $\delta_2 \in \F_{p^{2 s}}$ \\
\hline
\end{tabular}

\bigskip

\subsubsection{Case: $[G',G] \supsetneq Fratt(G')=\{e\}$.}

\begin{proposition} \label{casp2nontrivial}
Let $(C,G)$ be a big action satisfying $\mathcal{G}_*^{p^2}$  such that $[G',G] \neq\{e\}$.
We keep the notations introduced in section 5.1.
\begin{enumerate} 
\item
\begin{enumerate}
\item Then, $G=A_{\infty,1}$ is the unique $p$-Sylow subgroup of $A$.
\item For all $i$ in $\{1,2\}$, $f_i \in \Sigma_{i+1}-\Sigma_i$ and $m_i=1+i\,p^s$, with $p \geq 3$ and $s\in \{1,2\}$. \item Moreover, $v=s+1$.
 More precisely, $y \in V$ if and only if $\ell_{1,2}(y)^p-\ell_{1,2}(y)=0$ .
\end{enumerate}
\item There exists a coordinate $X$ for the projective line $C/G_2$ such that the functions $f_i$'s
are parametrized as follows:
 \end{enumerate}

\noindent \underline{(a) If $s=1$,}

\smallskip

\noindent \begin{tabular}{|c|c|c|}
\hline
& $p>3$ & $p=3$ \\
\hline
 $f_1$  &  $f_1(X)=X^{1+p}+a_2\, X^2$ & $f_1(X)=X^4+a_2 \, X$\\
\hline
$V$ & $V= Z(Ad_{f_1}) =Z(X^{p^2}+2\, a_2^p \,X^p +X)$&$V= Z(Ad_{f_1}) =Z(X^{9}+2\, a_2^3 \,X^3 +X)$\\
\hline
$f_2 $ & $f_2(X)= b_{1+2\,p} \, X^{1+2\,p} + b_{2+p} \, X^{2+p}  + b_{3} \,X^3+  b_{1}\,X$&
$f_2(X)= b_{7} \, X^{7} + b_{5} \, X^{5} + b_{1}\,X$\\
\hline
$b_{1+2\,p}$ & $ b_{1+2\,p} \in k^{\times}$ & $b_{7}^{16}=1$ \\
\hline
$a_2$ & $2\, a_2^p= -b_{1+2\,p}^{-p} (b_{1+2\,p}^{p^2}+b_{1+2\,p}) \Leftrightarrow b_{1+2\,p} \in V$& $2\, a_2^3= -b_{7}^{6} -b_{7}^{-2}$\\
\hline
$b_{2+p}$& $ b_{2+p}=-b_{1+2\,p}^p$& $ b_{5}=-b_{7}^3$\\
\hline
$b_{3}$ & $  3 \,b_{3}^p= b_{1+2\,p}^{-p} \, (b_{1+2\,p}^{2\,p^2}-b_{1+2\,p}^{2})$& \\
\hline
$b_{1}$ & $  b_{1} \in k$& $  b_{1} \in k$\\
\hline
$\ell_{1,2}$ & $\ell_{1,2}(y)=2\, (b_{1+2\,p} \, y^p-b_{1+2\,p}^p \, y)$
& $ \ell_{1,2}(y)=2\, (b_{7} \, y^3-b_{7}^3 \, y)$\\
\hline
\end{tabular} 

\bigskip

\noindent Therefore, for $p>3$, the solutions are parametrized by 2 algebraically independent variables over $\F_p$, namely $b_{1+2\,p} \in k^{\times}$ and $b_{1} \in k$. For $p=3$, as the monomial $X^3$ can be reduced mod $\wp(k[X])$, the parameter $b_{1+2\,p}$ satisfies an additional algebraic relation: $b_{7}^{16}=1$. Then, $b_{7}$ takes a finite number of values. 

\bigskip

\noindent \underline{In both cases ($p=3$ or $p>3$), }

\medskip

$$\frac{|G|}{g}= \frac{2\,p}{p-1} \,\frac{p^2}{1+2\,p}  \qquad \mbox{and} \qquad \frac{|G|}{g^2}=\frac{4}{(p^2-1)^2} \, \frac{p^2\, (p+1)^2}{(1+2\,p)^2}$$

\bigskip

\noindent  \underline{(b) If $s=2$ and $p>3$,}

\medskip

\noindent \begin{tabular}{|c|c|}
\hline
 $f_1$  &  $f_1(X)=X^{1+p^2}+a_{1+p}\, X^{1+p}+a_2\,X^2 $\\
 \hline
$Ad_{f_1}$ & $ X^{p^4} + a_{1+p}^{p^2} \, X^{p^3}+2\, a_2^{p^2} \, X^{p^2} +a_{1+p}^p \, X^p+X$\\
 \hline
$f_2$ &  $f_2(X)= b_{1+2\,p^2}\,  X^{1+2\,p^2} + b_{1+p+p^2} \, X^{1+p+p^2}  + b_{2+p^2} \, X ^{2+p^2} + b_{1+p^2}\, X ^{1+p^2} + b_{1+2\,p} \, X^{1+2\,p} $\\
& $+ b_{2+p} X^{2+p} + b_{1+p} X^{1+p} + b_{3} \, X^3+ b_{2} \, X^2 + b_{1}\,  X$ \\
\hline
$b_{1+2\,p}$&  $b_{1+2\,p} \in k^{\times}$\\
\hline
$b_{2+p^2}$& $b_{2+p^2} \in k^{\times}$\\
\hline
$b_{1+p+p^2}$ & $ b_{1+p+p^2} ^p = -2\,  b_{1+2\,p}^p \, ( b_{2+p^2}^p \, b_{1+2\,p}^{-p^2}+ b_{2+p^2}^{p-1}) $ \\
\hline
$\ell_{1,2}$ & $\forall \, y \in V, \quad \ell_{1,2}(y)=2\, b_{1+2\,p} \, y^{p^2} +b_{1+p+p^2} \, y^p+2\, b_{2+p^2} \, y$\\
\hline
$V$ & $V$ is an index p-subgroup of $Z(Ad_{f_1})$\\
&  $V=Z(2\, b_{1+2\,p}^p \, X^{p^3} +(b_{1+p+p^2}^p-2\, b_{1+2\,p}) \, X^{p^2} +(2 b_{2+p^2}^p- b_{1+p+p^2}) X^p -2b_{2+p^2} X)$\\
 \hline
$a_{1+p}$ & $a_{1+p}^{p^2} =  - b_{1+2\,p}^{p-p^2} - b_{1+2\,p}^p\, b_{2+p^2}^{-1} - b_{2+p^2}^{p^2}\, b_{1+2\,p}^{-p^3} - b_{2+p^2}^{p^2-p}$\\
\hline
$a_2$ & $ 2 \, a_2^{p^2}=  b_{2+p^2}^{p^2}\, b_{1+2\,p}^{-p^2} + b_{1+2\,p}\, b_{2+p^2}^{-1} + b_{2+p^2}^p\, b_{1+2\,p}^{p-2\,p^2} +2 \, b_{2+p^2}^{p-1}\, b_{1+2\,p}^{p-p^2} + b_{1+2\,p}^p b_{2+p^2}^{p-2}   $\\
\hline
$b_{1+2\,p}$ & $ b_{1+2\,p}^{p^2}= -b_{1+2\,p}^{2\,p-p^2}-b_{1+2\,p}^{2\,p}\, b_{2+p^2}^{-1}+
b_{2+p^2}^{2\,p^2}\, b_{1+2\,p}^{p^2-2\,p^3}+2\, b_{2+p^2}^{2\,p^2-p}\, b_{1+2\,p}^{p^2-p^3}+b_{1+2\,p}^{p^2} \, b_{2+p^2}^{2\,p^2-2\,p}$\\
\hline
$ b_{2+p}$ & $b_{2+p}^{p^2}= b_{2+p^2}^p\, b_{1+2\,p}^{2\,p-2\,p^2}+2\, b_{2+p^2}^{p-1}\, b_{1+2\,p}^{2\,p-p^2} +b_{1+2\,p}^{2\,p} \, b_{2+p^2}^{p-2} -b_{2+p^2}^{2\,p^2}\, b_{1+2\,p}^{-p^3} -b_{2+p^2}^{2\,p^2-p}$\\
\hline
$ b_{3} $ & $3 \, b_{3}^{p^2}=b_{2+p^2}^{2\,p^2}\, b_{1+2\,p}^{-p^2} - b_{2+p^2}^{2\,p}\, b_{1+2\,p}^{2\,p-3p^2} -3\, b_{2+p^2}^{2\,p-2}\, b_{1+2\,p}^{2\,p-p^2} -3 \, b_{2+p^2}^{2\,p-1}\, b_{1+2\,p}^{2\,p-2\,p^2}-b_{1+2\,p}^{2\,p}b_{2+p^2}^{2\,p-3} + b_{1+2\,p}^2\,b_{2+p^2}^{-1}$\\
\hline
$b_{1+p^2}$ & $ b_{1+p^2} \in Z(b_{2+p^2}^{p^2}\, b_{1+2\,p}^{-p^3} \, X^{p^3} - ( b_{2+p^2}^{p^2}\, b_{1+2\,p}^{-p^3} +b_{1+2\,p}^{p-p^2} +b_{2+p^2}^{p^2-p}) \, X^{p^2} +$\\
& $ (b_{1+2\,p}^{p-p^2} + b_{1+2\,p}^p\, b_{2+p^2}^{-1} + b_{2+p^2}^{p^2-p})\, X^p - b_{1+2\,p}^p\, b_{2+p^2}^{-1}\, X)$ \\
\hline 
$b_{1+p}$ & $b_{1+p}^{p^2}=  -(b_{1+2\,p}^{p-p^2} + b_{2+p^2}^{p^2}\,b_{1+2\,p}^{-p^3} + b_{2+p^2}^{p^2-p} )\, b_{1+p^2} ^{p^2} -b_{1+2\,p}^p\, b_{2+p^2}^{-1} \, b_{1+p^2}$\\
\hline
 $b_{2}$ & $2 \, b_{2}^{p^2} = (b_{2+p^2}^p\, b_{1+2\,p}^{p-2\,p^2} + b_{2+p^2}^{p-1}\,b_{1+2\,p}^{p-p^2-p} + b_{2+p^2}^{p^2}\, b_{1+2\,p}^{-p^2})
 b_{1+p^2}^{p^2} + $\\
& $( b_{1+2\,p}^p b_{2+p^2}^{p-2} + b_{2+p^2}^{p-1}\, b_{1+2\,p}^{p-p^2} + b_{1+2\,p}\, b_{2+p^2}^{-1}) \,b_{1+p^2}$ \\
\hline
$b_{1}$ & $b_{1} \in k$ \\
\hline
\end{tabular}

\bigskip

Therefore, for $p>3$, the solutions can be parametrized by 3 algebraically independent variables over $\F_p$, namely $b_{1+2\,p^2}\in k^{\times}$, $b_{2+p^2} \in k^{\times}$ and $b_{1} \in k$. One also finds a fourth parameter $b_{1+p^2}$ which runs  over an $\F_p$-subvector space of $k$, namely the set of zeroes of an additive separable polynomial whose coefficients are rational functions in $b_{1+2\,p^2}$ and $b_{2+p^2}$. So, for given $b_{1+2\,p}$ and $b_{2+p^2}$, the parameter $b_{1+p^2}$ takes a finite number of values.\\
\medskip 

\noindent \underline{For $p=3$,} 
\medskip

\noindent 
$$
\begin{array}{|ll|}
\hline
f_1(X)=&X^{10}+a_4\, X^{4}+a_2\,X^2  \\
\hline
f_2(X)=& b_{19} \,X^{19} + b_{13} \,X^{13}  + b_{11} \,X ^{11} + b_{10}\, X ^{10} + b_{7}\, X^{7} + b_{5} \, X^{5}+\\
 &b_{4} \,X^{4} + b_{2} \,X^2 + b_{1} \,X \\
\hline
\end{array}
$$

\medskip

\noindent with $a_4$, $a_2$, $b_{13}$, $b_{7}$, $b_{5}$, $b_{3}$ and $b_{2}$ satisfying the same relations as above.
But, this time, the parameters $b_{19}$ and $b_{11}$ are linked through an algebraic relation, namely:

$$b_{11}^{18}\, b_{19}^{-9} - b_{11}^{6}\, b_{19}^{-21} -b_{19}^{6}\, b_{11}^{3} + b_{19}^2\, b_{11}^{-1}=0$$

\medskip

\noindent \underline{In both cases ($p=3$ or $p>3$), }

$$ \frac{|G|}{g}=\frac{2\,p}{p-1} \,\, \frac{p^2}{1+2\,p} \qquad \mbox{and} \qquad \frac{|G|}{g^2}=\frac{4}{(p^2-1)^2} \, \, \frac{p\, (p+1)^2}{(1+2\,p)^2}$$

\end{proposition}

\begin{remarque}
One can now answer the second problem raised in \cite{MR2} (section 6). Indeed, one notices that the family obtained for $s=2$ is larger than the one obtained after the additive base change: $X=Z^p+c\,Z$, $c\in k-\{ 0\}$ (see \cite{MR} Prop. 3.1) applied to the case $s=1$. Indeed, such a base change does not produce any monomial $Z^{1+p^2}$ in $f_2(Z)$.
\end{remarque}

\textbf{A few special cases.}
\begin{enumerate}
\item When $s=1$ and $p>3$, the special case $a_2=0$  corresponds to the
parametrization of the extension $K_S^m/K$ given by Auer (cf . \cite{Au} Prop. 8.9 or \cite{MR} section 6), namely
\begin{description}
\item $f_1(X)= a\, X^{1+p}$ with $a^p+a=0$, $a\neq 0$.
\item $f_2(X)=a^2 \, X^{2\,p} \, (X-X^{p^2})$.
\end{description}
\item When $s=2$, the special case $b_{1+p^2} \in \F_p$ leads to $b_{1+p}= b_{1+p^2} \, a_{1+p}$ and $b_{2}=b_{1+p^2} \, a_2$. So, one can replace $f_2$ by  $f_2(X)-b_{1+p^2} \, f_1(X)$, which eliminates the monomials $X^{1+p^2}$, $X^{1+p}$ and $X^2$.
\end{enumerate}

\noindent \textbf{Proof of Proposition \ref{casp2nontrivial}:}
\begin{enumerate}
\item
As $\ell_{1,2} \neq 0$, the group $G$ satisfies the third condition of \cite{MR2} (Prop. 5.2). 
Then, the equality $G=A_{\infty,1}$ derives from \cite{MR2} (Cor. 5.7). The unicity of the $p$-Sylow subgroup is explained in Remark \ref{italien}. The second and third assertions come from \cite{MR2} (Thm. 5.6). Moreover, the description of $V$ displayed in $(c)$ is due to\cite{MR2} (Prop 2.9.2). 
It remains to show that $s=1$ or $s=2$. Using formula \eqref{genus}, we compute $g=\frac{(p-1)}{2}\,(p^{s}+p\, (m_2-1))=\frac{(p-1)}{2}\,p^s (1+2 \,p).$
 As $|G|=p^{3+s}$, condition $(*)$ requires: $\frac{4}{(p^2-1)^2} \leq  \frac{|G|}{g^2}  =\frac{4}{(p^2-1)^2}\,  \frac{(p+1)^2}{p^{s-3} \, (1+2\,p)^2}.$ It follows that $3-s>0$, i.e. $1 \leq s \leq 2$. 

\item We merely explain the case $s=1$. One can find a coordinate $X$ of the projective line $C/G_2$ such that $f_1(X)= X \,S_1(X)=X\, (X^p + a_2\,X)$ (cf. \cite{MR2} Cor. 2.12). 
Then, $Ad_{f_1}=F^2+2\, a_2^p \, F + I$ (cf. \cite{MR2} Prop. 2.13). As $V \subset Z(Ad_{f_1})$ and $dim_{\F_p} \, Z(Ad_{f_1})=2=s+1=v$, we deduce that
$V = Z(Ad_{f_1})$. As $f_2 \in \Sigma_3-\Sigma_2$ with $deg \, f_2=1+2\,p^s$
and as the functions $f_i$'s are supposed to be reduced mod $\wp(k[X])$, equation \eqref{toto} reads:

$$
 \begin{array}{lll}
 \forall \, y \in V, &\quad f_2(X+y)-f_2(X)= \ell_{1,2}(y)\, f_1(X) \quad \mod \, \wp(k[X]) &\\ 
& &\\
\mbox{with} & f_1(X)= X^{1+p}+ a_2 \, X^2&\\
  \mbox{and}&
  f_2(X)= b_{1+2\,p} \, X^{1+2\,p} + b_{2+p} \,X^{2+p}+ b_{1+p} \,X^{1+p} + b_{3} \,X^3 + b_{2}
\, X^2+b_{1} \, X & \mbox{for} \, \,  p>3\\
    (\mbox{resp.}&f_2(X)= b_{1+2\,p} \, X^{1+2\,p} + b_{2+p} \,X^{2+p}+ b_{1+p} \,X^{1+p} + b_{2} \,X^2 +b_{1}
    \,X& \mbox{for} \, \,  p=3)\\
\end{array}
$$

Then, calculation gives the relations gathered in the table. In particular, we find: $f_2(X)= b_{1+2\,p} X^{1+2\,p}+ b_{2+p} \,   X^{2+p} +b_{3}\, X^3 +  b_{1}\, X + b_{1+p}\,  f_1(X)$ with  $b_{1+p} \in \F_p$. Since we are working in the $\F_p$-space generated by $f_1(X)$ and $f_2(X)$, we can replace $f_2(X)$ with $f_2(X)-b_{1+p} \, f_1(X)$, hence the expected formula.
We solve the case $s=2$ in the same way. $\square$
\end{enumerate}

\subsection{Third case: big actions satisfying $\mathcal{G}_*^{p^3}$.}

\subsubsection{Preliminaries.}
The idea is to use, as often as possible, the results obtained in the preceding section.

\begin{remarque} \label{deb}
Let $(C,G)$ be a big action with $G'(=G_2) \simeq (\Z/p\Z)^3$.
We keep the notations introduced in section 5.1.
\begin{enumerate}
\item Let $C_{1,2}$ be the curve parametrized by the two equations: $W_i^p-W_i=f_i(X)$, with $i \in \{1,2\}$,  and let $K_{1,2}:=k(C_{1,2})$ be the function field of this curve. Then, $K_{1,2}/k(X)$ is a Galois extension with group $\Gamma_{1,2} \simeq (\Z/p\Z)^2$. Moreover, the group of translations by $V$: $\{X\rightarrow X+y, y \in V\}$ extends to an automorphism $p$-group of $C_{1,2}$ say $G_{1,2}$, with the following exact sequence:
$$0 \longrightarrow \Gamma_{1,2} \longrightarrow G_{1,2} \longrightarrow V \longrightarrow 0$$
Let $A_{1,2}$ be the $\F_p$-subvector space of $A$ generated by the classes of $f_1(X)$ and $f_2(X)$. Let $H_{1,2} \subsetneq G_2$ be the orthogonal of $A_{1,2}$ with respect to the Artin-Schreier pairing. Then, 
$C_{1,2}=C/H_{1,2}$ and $G_{1,2}=G/H_{1,2}$. Furthermore, as $A_{1,2}$ is stable under the action of $\rho$,
its dual $H_{1,2}$ is stable by the dual representation $\phi$, i.e. by conjugation by the elements of $G$ (see section 5.1). It follows that $H_{1,2} \subsetneq G_2$ is a normal subgroup in $G$. So, by \cite{MR} (Lemma 2.4), the pair $(C_{1,2}, G_{1,2})$ is a big action with second ramification group isomorphic to $(\Z/p\Z)^2$. 
\item Likewise, if $\ell_{2,3}=0$, the $\F_p$-subvector space of $A$ generated by the classes of ${f_1(X)}$ and $f_3(X)$  is also stable by $\rho$ (see matrix $L(y)$ in section 5.1). So, the two equations: $W_i^p-W_i=f_i(X)$, with $i \in \{1,3\}$, also
parametrize a big action, say $(C_{1,3}, G_{1,3})$, with second ramification group isomorphic to $(\Z/p\Z)^2$. 
\item Similarly, if $\ell_{1,2}=\ell_{1,3}=0$, the $\F_p$-subvector space of $A$ generated by the classes of $f_2(X)$ and $f_3(X)$ is stable by $\rho$ (see matrix $L(y)$ in section 5.1). So, the two equations: $W_i^p-W_i=f_i(X)$, with $i \in \{2,3\}$, also
parametrize a big action, say $(C_{2,3}, G_{2,3})$, with second ramification group isomorphic to $(\Z/p\Z)^2$. 
\end{enumerate}
\end{remarque}

\begin{lemme} \label{G12}
Let $(C,G)$ be a big action satisfying $\mathcal{G}_*^{p^3}$. Let $(C_{1,2}, G_{1,2})$ be defined as in Remark \ref{deb}.
We keep the notations introduced in section 5.1.
\begin{enumerate}
\item Then, $(C_{1,2}, G_{1,2})$ is a big action satisfying  $\mathcal{G}_*^{p^2}$. 
\item If $\ell_{1,2}=0$, then  $m_1=m_2=1+ p^s$, with $s\geq 2$.
\item If $\ell_{1,2}\neq 0$, then $m_1=1+p^s$, $m_2=1+2\, p^s$, with $s \in \{1,2\}$ and $p \geq 3$. In this case, $v=s+1$.
\end{enumerate}
\end{lemme}

\noindent \textbf{Proof:} 
\begin{enumerate}
\item Use Remark \ref{deb} and \cite{LM} (Prop. 8.5 (ii)) to see that condition $(*)$ is still satisfied.
\item We  deduce from Proposition \ref{casp2trivial} that $m_1=1+p^{s_1}$ and $m_2=1+p^{s_2}$ with $s_2=s_1$ or $s_2=s_1+1$. Assume that $s_2=s_1+1$. Then, $m_3\geq m_2=1+p^{s_1+1}$. We compute the genus by means of \eqref{genus}:
$g=\frac{p-1}{2} \, (p^{s_1}+p^{1+s_2}+p^2\, (m_3-1)) \geq \frac{p-1}{2} \, p^{s_1} \, (1+p^2+p^3)$. 
Besides, by \cite{MR} (Thm. 2.6),  $V \subset Z(Ad_{f_1})$, so $|G|=p^{3+v} \leq p^{3+2s_1}$. Thus, 
$\frac{|G|}{g^2} \leq \frac{4}{(p^2-1)^2} \, \, \frac{p^3\, (p+1)^2}{(1+p^2+p^3)^2}< \frac{4}{(p^2-1)^2} $,
which contradicts condition $(*)$. It follows that $s_2=s_1\geq 2$. 
\item Apply Proposition \ref{casp2nontrivial} to $(C_{12}, G_{12})$.
$\square$
\end{enumerate}

\begin{remarque} \label{G23}
Let $(C,G)$ be a big action satisfying $\mathcal{G}_*^{p^3}$. Assume that $\ell_{1,2}=\ell_{1,3}=0$. Then, the results of Lemma \ref{G12} also hold for the big action $(C_{2,3}, G_{2,3})$ as defined in Remark \ref{deb}.
\end{remarque}

\begin{lemme} \label{G13}
Let $(C,G)$ be a big action satisfying $\mathcal{G}_*^{p^3}$. We keep the notations introduced in section 5.1.
Assume that $\ell_{2,3}=0$.
 Let $(C_{1,3}, G_{1,3})$ be defined as in Remark \ref{deb}.
\begin{enumerate}
\item Then, $(C_{1,3}, G_{1,3})$ is a big action satisfying  $\mathcal{G}_*^{p^2}$. 
\item If $\ell_{1,3}=0$, then $\ell_{1,2}=0$ and  $m_1=m_2=m_3=1+ p^s$ with $s\geq 2$. In this case, $v=2\,s$.
\item If $\ell_{1,3}\neq 0$, then $m_1=1+p^s$, $m_3=1+2\, p^s$, with $s \in \{1,2\}$ and $p \geq 3$. In this case, $v=s+1$.
\end{enumerate}
\end{lemme}

\noindent \textbf{Proof:} 
\begin{enumerate}
\item Use Remark \ref{deb} and \cite{LM} (Prop. 8.5 (ii)).
\item As $\ell_{1,3}=0$, we deduce from Proposition \ref{casp2trivial} that $m_1=1+p^{s}$ and $m_3=1+p^{s_3}$ with $s_3=s$ or $s_3=s+1$.
\begin{enumerate}
\item \emph{We show that $\ell_{1,2}=0$}.\\
Assume that $\ell_{1,2} \neq 0$. Then, Lemma \ref{G12} applied to $(C_{1,2},G_{1,2})$ implies
 $m_2=1+2\, p^{s}$ with $s \in \{1,2\}$ and $p \geq 3$. Moreover, $v=s+1$. As $m_2 \leq m_3$, there are two possibilities:
\begin{enumerate}
\item \emph{$s=1$ and $s_3=s+1=2$.}, i.e.  $m_1=1+p$, $m_2=1+2\,p$, $m_3=1+p^2$ and $v=2$.
Then, $\frac{|G|}{g^2}=\frac{4}{(p^2-1)^2} \, \, \frac{p^3\, (p+1)^2}{(1+2\,p+p^3)^2} <\frac{4}{(p^2-1)^2}$, which contradicts condition $(*)$.
\item \emph{$s=2$ and $s_3=s+1=3$.} i.e.  $m_1=1+p^2$, $m_2=1+2\,p^2$, $m_3=1+p^3$ and  $v=3$.
Then, $\frac{|G|}{g^2}=\frac{4}{(p^2-1)^2} \, \, \frac{p^2\, (p+1)^2}{(1+2\,p+p^3)^2} <\frac{4}{(p^2-1)^2}$, which also contradicts condition $(*)$.
\end{enumerate}
As a consequence, $\ell_{1,2}=0$. 
\item \emph{We deduce that $m_1=m_2=1+p^{s}$ with $s\geq 2$.}\\
Lemma \ref{G12} applied to $(C_{1,2},G_{1,2})$ implies
$m_1=m_2=1+p^{s}$ with $s\geq 2$. In particular, $g =\frac{p-1}{2} \, p^s \, (1+p+p^{2+s_3-s})$ and $\frac{|G|}{g^2}=\frac{4}{(p^2-1)^2} \, \frac{p^{3+v-2s} \, (p+1)^2}{(1+p+p^{2+s_3-s})^2} $. 
\item \emph{We show that $v=2\,s_3$ and conclude that $s_3=s$.} \\
Assume that $v \leq 2\, s_3-3$.
Then,
$\frac{|G|}{g^2} < \frac{4}{(p^2-1)^2} \, \frac{p^{2s_3-2s} \, (p+1)^2}{p^{4+2s_3-2s}}< \frac{4}{(p^2-1)^2}$
which contradicts condition $(*)$. Therefore, $2\, s_3-2 \leq v \leq 2\,s \leq 2\, s_3$.
Assume that $v \leq 2\, s_3-3$.
Then,
$\frac{|G|}{g^2} < \frac{4}{(p^2-1)^2} \, \frac{p^{2s_3-2s} \, (p+1)^2}{p^{4+2s_3-2s}}< \frac{4}{(p^2-1)^2}$
which contradicts condition $(*)$.
Assume that $v=2\,s_3-1$. So, $v$ is odd and $2\,s_3-2 <v \leq 2s \leq 2\,s_3$ implies $s_3=s$ and $v=2s-1$. Then, $\frac{|G|}{g^2}= 
\frac{4}{(p^2-1)^2} \, \frac{p^{2} \, (p+1)^2}{(1+p+p^{2})^2}< \frac{4}{(p^2-1)^2}$, which is excluded.
Now, assume that $v=2\,s_3-2$. Then, $2\,s_3-2 =v \leq 2s \leq 2\,s_3$ implies $s_3=s$ or $s_3=s+1$.
In the first case, $v=2s-2$ and
$\frac{|G|}{g^2}=
 \frac{4}{(p^2-1)^2} \, \frac{p \, (p+1)^2}{(1+p+p^{2})^2}< \frac{4}{(p^2-1)^2}$.
In the second case, $v=2s$ and 
$\frac{|G|}{g^2}= \frac{4}{(p^2-1)^2} \, \frac{p^3 \, (p+1)^2}{(1+p+p^{3})^2} < \frac{4}{(p^2-1)^2}$.
In both cases, we obtain a contradiction.
We gather that $v=2\,s_3$. Applying \cite{MR2} (Prop. 4.2), we conclude that $s=s_3$. 
\end{enumerate}
\item Apply Proposition \ref{casp2nontrivial} to $(C_{13}, G_{13})$. $\square$
\end{enumerate}

\subsubsection{Case: $[G',G]=Fratt(G')=\{e\}$.}

\begin{proposition} \label{p3cas1}
Let $(C,G)$ be a big action satisfying $\mathcal{G}_*^{p^3}$ such that $[G',G]=\{e\}$. 
We keep the notations introduced in section 5.1.
\begin{enumerate}
\item Then, $G=A_{\infty,1}$ is a special group of exponent $p$ (resp. $p^2$) for $p>2$ (resp. $p=2$) and order $p^{3+2\, s_1}$.
More precisely, $G$ is a central extension of its center $Z(G)=G'$ by the elementary abelian $p$-group $V=Z(Ad_{f_1})= Z(Ad_{f_2})=Z(Ad_{f_3})$:
$$ 0 \longrightarrow Z(G)=G'\simeq (\Z/p \Z)^3 \longrightarrow G \stackrel{\pi}{\longrightarrow} Z(Ad_{f_1})= Z(Ad_{f_2})=Z(Ad_{f_3}) \simeq (\Z / p \Z) ^{2s_1} \longrightarrow 0$$ 
Furthermore, $G$ is a $p$-Sylow subgroup of $A$, which is normal except when $C$ is birational to the Hermitian curve: $W^q-W=X^{1+q}$, with $q=p^3$. 
\item There exists a coordinate $X$ for the projective line $C/G_2$, $s \geq 2$, $d \geq 2$ dividing $s$, and $\gamma_2$, $\gamma_3$ in $\F_{p^d}-\F_p$ linearly independent over $\F_p$, $b_1 \in k$, $c_1 \in k$ such that: \\

\medskip
\begin{tabular}{|c|c c c|}
\hline
$f_1$ & $f_1(X)=X\, S_1(X)$ &\mbox{with} & $S_1(F)=\sum_{j=0}^{s/d}  \, a_{jd} \, F^{jd} \in k\{F\} \quad a_s=1$ \\
\hline
$f_2$ & $f_2(X)=X\, S_2(X)+b_1\, X$ &\mbox{with} & $S_2= \gamma_2 \, S_1   $\\
\hline
$f_3$ & $f_3(X)=X\, S_3(X)+c_1\, X$ &\mbox{with} & $S_3= \gamma_3 \, S_1 $\\
\hline
$V$ & $V=Z(Ad_{f_1})= Z(Ad_{f_2})=Z(Ad_{f_3})$ & & \\
\hline
\end{tabular}

\medskip

Therefore, the solutions can be parametrized by $s+4$ algebraically independent variables over $\F_p$, namely the $s$ coefficients of $S$, $\gamma_2 \in \F_{p^d}-\F_p$, $\gamma_3 \in \F_{p^d}-\F_p$ , $b_1 \in k$ and $c_1 \in k$.

\bigskip
\noindent Moreover,
$$ \frac{|G|}{g}=\frac{2\,p}{p-1} \, \,  \frac{p^{s}}{1+p+p^2}\quad \mbox{and} \quad 
\frac{|G|}{g^2}=\frac{4}{(p^2-1)^2} \, \,\frac{p^3(p+1)^2}{(1+p+p^2)^2} $$
\end{enumerate}
\end{proposition}

\noindent \textbf{Proof:}
As $\ell_{1,2}=\ell_{2,3}=\ell_{1,3}=0$, the second point of Lemma \ref{G13} first implies $v=2\,s_3$.
Applying \cite{MR2} (Prop. 4.2), we gather that $s_1=s_2=s_3$,
that $V=Z(Ad_{f_1})= Z(Ad_{f_2})=Z(Ad_{f_3})$ and we get the expected formulas for the functions $f_i's$.
Moreover, it follows from \cite{MR2} (Prop. 4.3 and Rem. 4.5) that $G=A_{\infty,1}$ is a special group.
The unicity of the $p$-Sylow subgroup is discussed in Remark \ref{italien}.
$\square$

\subsubsection{Case: $[G',G] \supsetneq Fratt(G')=\{e\}$.}

\begin{lemme}
Let $(C,G)$ be a big action satisfying $\mathcal{G}_*^{p^3}$ such that $[G',G] \neq \{e\}$.
We keep the notations introduced in section 5.1.
Then, one cannot have $\ell_{1,2}=\ell_{2,3}= 0$. 
\end{lemme}

\noindent \textbf{Proof:}
Assume that $\ell_{1,2}=0$ and $\ell_{2,3}=0$. Since the representation $\rho$ is non trivial, $\ell_{1,3}\neq 0$.
The second point of Lemma \ref{G12} shows that $m_1=m_2=1+p^s$ with $s\geq 2$. The third point of Lemma \ref{G13}  implies that $m_3=1+2\, p^{s}$ with $p \geq 3$ and $s\in \{1,2\}$. Moreover, $v=s+1$. As $s\geq 2$, we obtain:
 $\frac{|G|}{g^2} =\frac{4}{(p^2-1)^2} \, \frac{(p+1)^2 \, p^2}{(1+p+2\,p^2)^2} <\frac{4}{(p^2-1)^2} $, hence
a contradiction. As a conclusion, either $\ell_{1,2}\neq 0$ or $\ell_{2,3}\neq 0$. $\square$
\medskip

As a consequence, there are 3 cases to study: 
\begin{description}
\item $\ell_{1,2}\neq 0$ and $\ell_{2,3}= 0$ (cf. Proposition \ref{p3cas1}).
\item $\ell_{1,2}=0$ or $\ell_{2,3}\neq  0$ (cf. Proposition \ref{p3cas2}).
\item  $\ell_{1,2}\neq 0$ or $\ell_{2,3}\neq 0$ (cf. Proposition \ref{p3cas3}).
\end{description}

\begin{proposition} \label{p3cas1}
Let $(C,G)$ be a big action satisfying $\mathcal{G}_*^{p^3}$ such that $[G',G] \neq \{e\}$.
We keep the notations introduced in section 5.1. Assume that $\ell_{1,2}\neq 0$ and $\ell_{2,3}= 0$. 
\begin{enumerate}
\item Then, $p\geq 5$ and there exists a coordinate $X$ for the projective line $C/G_2$ such that
the functions $f_i$'s can be parametrized as follows:
\medskip

\noindent \begin{tabular}{|c|c|}
\hline
 $f_1$  &  $f_1(X)=X^{1+p}+a_2\, X^2$ \\
\hline 
$V$ & $V= Z(Ad_{f_1}) =Z(X^{p^2}+2\, a_2^p \,X^p +X)$\\
\hline
\hline
$f_2 $ & $f_2(X)= b_{1+2\,p} \, X^{1+2\,p} + b_{2+p} \, X^{2+p}  + b_{3} \,X^3+  b_{1}\,X$\\
\hline
$b_{1+2\,p}$ & $ b_{1+2\,p} \in k^{\times}$ \\
\hline
$a_2$ & $2\, a_2^p= -b_{1+2\,p}^{-p} \, (b_{1+2\,p}+b_{1+2\,p}^{p^2}) \Leftrightarrow b_{1+2\,p} \in V$\\
\hline 
$V$&$V=Z(X^{p^2} -b_{1+2\,p}^{-p} \, (b_{1+2\,p}+b_{1+2\,p}^{p^2}) \,X^p +X)$\\
\hline
$b_{2+p}$& $ b_{2+p}=-b_{1+2\,p}^p$\\
\hline
$b_{3}$ & $  3 \,b_{3}^p= b_{1+2\,p}^{-p} \, (b_{1+2\,p}^{2\,p^2}-b_{1+2\,p}^{2})$ \\
\hline
$b_{1}$ & $  b_{1} \in k$\\
\hline
$\ell_{1,2}$ & $\ell_{1,2}(y)=2\, (b_{1+2\,p} \, y^p-b_{1+2\,p}^p \, y)$\\
\hline
\hline
$f_3 $ & $f_3(X)= c_{1+2\,p} \, X^{1+2\,p} + c_{2+p} \, X^{2+p}  + c_3 \,X^3+  c_1\,X$\\
\hline
$c_{1+2\,p}$ & $ c_{1+2\,p} \in k^{\times}$ \\
\hline
$c_{1+2\,p} $ & $ c_{1+2\,p} \in V$,  $c_{1+2\,p}$ \mbox{and} $b_{1+2\,p}$ $\F_p$-\mbox{independent}  \\
\hline
$c_{2+p}$& $ c_{2+p}=-c_{1+2\,p}^p$\\
\hline
$c_3$ & $  3 \,c_3^p= -c_{1+2\,p}^{-p} \, (c_{1+2\,p}^{2\,p^2}+c_{1+2\,p}^{2})$ \\
\hline
$c_1$ & $  c_1 \in k$\\
\hline
$\ell_{1,3}$ & $\ell_{1,3}(y)=2\, (c_{1+2\,p} \, y^p-c_{1+2\,p}^p \, y)$\\
\hline
$\ell_{2,3}$ & $\ell_{2,3}(y)=0$\\
\hline

\end{tabular}

\bigskip

\noindent Therefore, the solutions are parametrized by 4 algebraically independent variables over $\F_p$,
namely $b_{1+2\,p} \in k^{\times}$, $c_{1+2\,p} \in k^{\times}$, $b_{1} \in k$ and $c_1 \in k$. \\
\medskip
\noindent Moreover, 
$$ \frac{|G|}{g}=\frac{2\,p}{p-1} \,\,   \frac{p^3}{1+2\,p+2\,p^2}\quad \mbox{and} \quad 
\frac{|G|}{g^2}= \frac{4}{(p^2-1)^2} \, \,\frac{p^3(p+1)^2}{(1+2\,p+2\,p^2)^2}$$
\item In this case, $G=A_{\infty,1}$ is the unique $p$-Sylow subgroup of $A$.

\end{enumerate}
\end{proposition}

\noindent \textbf{Proof:} 
\begin{enumerate}
\item Lemma \ref{G12} first shows that $m_1=1+p^s$, $m_2=1+2\,p^s$, 
with $p \geq 3$ and $s\in \{1,2\}$. Moreover, $v=s+1$. 
As  $\ell_{1,2}\neq 0$ and $\ell_{2,3}=0$, the second point of Lemma \ref{G13} imposes $\ell_{1,3} \neq 0$. Then, Lemma \ref{G13} shows that $m_3=1+2\,p^{s}$. If $s=2$, $m_1=1+p^2$, $m_2=m_3=1+2\,p^3$ and $v=3$. So $\frac{|G|}{g^2}=\frac{4}{(p^2-1)^2} \, \, \frac{p^2 \, (p+1)^2}{(1+2\,p+2\,p^2)^2} <\frac{4}{(p^2-1)^2}$, which contradicts condition $(*)$. It follows that $s=1$. In this case, $m_1=1+p$, $m_2=m_3=1+2\,p$,  $v=2$ and $\frac{|G|}{g^2}=\frac{4}{(p^2-1)^2} \, \, \frac{p^3 \, (p+1)^2}{(1+2\,p+2\,p^2)^2}$. Therefore, condition $(*)$ is satisfied as soon as $p\geq 5$.
The parametrization of the functions $f_i$'s then derives from Proposition \ref{casp2nontrivial}.
Furthermore, the third condition (cf. Recall 4.2.1-c) imposed on the degree of the functions $f_i$'s requires that the parameters $b_{1+2\,p}$ and $c_{1+2\,p}$ are linearly independent over $\F_p$. 
\item The equality $G=A_{\infty,1}$ derives from the maximality of $V=Z(Ad_{f_1})$ (see Proposition \ref{prop11}). The unicity of the $p$-Sylow subgroup is due to Remark \ref{italien}. $\square$
\end{enumerate}

\begin{proposition} \label{p3cas2}
Let $(C,G)$ be a big action satisfying $\mathcal{G}_*^{p^3}$ such that $[G',G] \neq \{e\}$.
We keep the notations introduced in section 5.1. 
Assume that $\ell_{1,2}= 0$ and $\ell_{2,3}\neq 0$. 
\begin{enumerate}
\item Then, $p\geq 5$ and there exists a coordinate $X$ for the projective line $C/G_2$ such that
the functions $f_i$'s can be parametrized as follows:
\medskip

\noindent \begin{tabular}{|c|c|}
\hline
 $f_1$  &  $f_1(X)= X^{1+p^2}+ a_2 \, X^2$\\
\hline
\hline
$f_2 $ & $f_2(X)= \gamma_2\, ( \, X^{1+p^2}+ a_2 \, X^2) + b_1 \, X$\\
\hline
$b_1$ & $b_1 \in k$ \\
\hline
 $\gamma_2$ &   $\gamma_2 \in \F_{p^2} - \F_p$ \\
\hline
$V$ & $V= Z(Ad_{f_1})= Z(Ad_{f_2}) = Z( X^{p^4} + 2\,  a_2 ^{p^2} X ^{p^2} + X)$\\
\hline
\end{tabular}

\bigskip

\noindent \underline{First case: $b_1 \neq 0$}\\

\noindent \begin{tabular}{|c|c|}
\hline
 $f_3$  &  $f_3(X)= c_{1+2\,p^2} \, X^{1+2\,p^2} + c_{2+p^2} \, X ^{2+p^2} + c_{1+p^2} \, X ^{1+p^2}$ \\
& $+ c_{1+p} \, X^{1+p} + c_{3} \, X^3+  c_{2} X^2 + c_{1} \,  X$\\
\hline
$c_{1+2\,p^2}$ & $ c_{1+2\,p^2} \in k^{\times}$\\
\hline
$a_2$ & $2\,a_2^{p^2}=-c_{1+2\,p^2}^{-p^2} \, (c_{1+2\,p^2}^{p^4}+c_{1+2\,p^2}) \Leftrightarrow c_{1+2\,p^2} \in V$\\
\hline 
$V$ &$  V= Z( X^{p^2} -c_{1+2\,p^2}^{-p} \, (c_{1+2\,p^2}+c_{1+2\,p^2}^{p^2})\, X^p +X)$ \\ 
\hline
$c_{2+p^2} $ & $ c_{2+p^2} = -c_{1+2\,p^2} ^{p^2}$\\
\hline
$c_{3}$ & $ 3\,  c_{3}^{p^2}= - c_{1+2\,p^2}^{p^2} \, (3\, c_{1+2\,p^2}^{2\, p^4} +4\, c_{1+2\,p^2}^{1+p^4} +c_{1+2\,p^2}^2) $\\
\hline
$e:=c_{1+p^2} - c_{1+p^2}^{p^2}$ & $ e \in Z\, ((c_{1+2\,p^2}^{p^7-p^3} +1+c_{1+2\,p^2}^{p-p^5} +c_{1+2\,p^2}^{p^7+p-p^5-p^3}) \, X^{1+p^4}$\\
&$- X^{1+p^2} - X^{p^2}-X -1)$\\
\hline
$b_1 $ & $ b_1^{p^5-p^4+p^3-p^2} = -e ^{p^3-1}$ \\
\hline
$c_{1+p}$ & $  c_{1+p}^{p+p^3} =-e^{1+p} $ \\
\hline
$c_{2} $ &$ 4\, c_{2}^{p^3\,(p-1)^2\, (p^2+1)} = 
c_{1+p^2}^{p^3\, (p^2+1)} \, $\\
&$+(c_{1+2\,p^2}^{p^7-p^3} +1+c_{1+2\,p^2}^{p-p^5} +c_{1+2\,p^2}^{p^7+p-p^5-p^3} )\, (c_{1+p^2} - c_{1+p^2}^{p^2})^{p^3+p^2+p-1}$ \\
\hline
$c_1$ & $c_1\in k$\\
\hline
\end{tabular}

\bigskip

\noindent Therefore, the solutions can be parametrized by 3 algebraically independent variables over $\F_p$, namely $c_{1+2\,p^2}\in k^{\times}$, $c_{1} \in k$ and $\gamma_2 \in \F_{p^2}-\F_p$. One also finds a fourth parameter $e:=c_{1+p^2}-c_{1+p^2}^{p^2}$ which runs  over the set of zeroes of a polynomial whose coefficients are rational functions in $c_{1+2\,p^2}$. So, for a given $c_{1+2\,p^2}$, the parameter $e$ takes a finite number of values.\\

\medskip

\noindent \underline{Second case: $b_1 = 0$}\\

\noindent \begin{tabular}{|c|c|}
\hline
 $f_3$  &  $f_3(X)= c_{1+2\,p^2}\,  X^{1+2\,p^2} + c_{2+p^2} \,  X ^{2+p^2} + c_{3}\,  X^3$\\
\hline
$c_{1+2\,p^2}$ & $ c_{1+2\,p^2} \in k^{\times}$\\
\hline
$a_2$ & $2\,a_2^{p^2}=-c_{1+2\,p^2}^{-p^2} \, (c_{1+2\,p^2}^{p^4}+c_{1+2\,p^2}) \Leftrightarrow c_{1+2\,p^2} \in V$\\
\hline
$V$ &$ V=Z(X^{p^2}-c_{1+2\,p^2}^{-p} \, (c_{1+2\,p^2}+c_{1+2\,p^2}^{p^2}) \,X^p +X)$ \\ 
\hline
$c_{2+p^2} $ & $ c_{2+p^2} = -c_{1+2\,p^2} ^{p^2}$\\
\hline
$c_{3}$ & $ 3\,  c_{3}^{p^2}= - c_{1+2\,p^2}^{p^2} \, (3\, c_{1+2\,p^2}^{2\, p^4} +4\, c_{1+2\,p^2}^{1+p^4} +c_{1+2\,p^2}^2) $\\
\hline
$c_1$ & $c_1\in k$\\
\hline
\end{tabular}

\bigskip 

\noindent In this case, the solutions can be parametrized by 3 algebraically independent variables over $\F_p$, namely $c_{1+2\,p^2}\in k^{\times}$, $c_1 \in k$ and $\gamma_2 \in \F_{p^2}-\F_p$.

\medskip

\noindent In both cases, 
$$ \frac{|G|}{g}=\frac{2\,p}{p-1} \, \,  \frac{p^4}{1+p+2\,p^2}\quad \mbox{and} \quad 
\frac{|G|}{g^2}= \frac{4}{(p^2-1)^2} \,\,  \frac{p^3\,(p+1)^2}{(1+p+2\,p^2)^2}$$
\item Moreover, $G=A_{\infty,1}$ is the unique $p$-Sylow subgroup of $A$.
\end{enumerate}

\end{proposition}
 
\noindent \textbf{Proof:} 
\begin{enumerate}
\item \begin{enumerate}
\item \emph{We describe $f_1$, $f_2$ and $V$.} \\
Lemma \ref{G12} first implies that $m_1=m_2=1+p^s$, with $s\geq 2$. More precisely, we deduce from Proposition \ref{casp2trivial} that $f_1(X)=X\, S_1(X)$ and $f_2(X)=\gamma_2\, X \, S_1(X)+b_1\, X$, where $S_1$ is a monic additive polynomial with degree $s$ in $F$, $b_1 \in k$ and $\gamma_2 \in \F_{p^d}-\F_p$ with $d$ an integer dividing $s$. Moreover, $v=2\,s$ and $V=Z(Ad_{f_1})=Z(Ad_{f_2})$. 
\item \emph{We show that $\ell_{1,3} \neq 0$.} \\
Indeed, assume that $\ell_{1,3}=0$. Then, we deduce from Remark \ref{G23} that $m_3=1+2\, p^s$, with $s \in \{2,3\}$ and $p \geq 3$. Moreover, $v=s+1$. As $s \neq 1$, it follows that $s=2$ and $\frac{|G|}{g^2}= \frac{4}{(p^2-1)^2} \, \frac{p^2\, (p+1)^2}{(1+p+2\, p^2)^2} < \frac{4}{(p^2-1)^2}$, which contradicts condition $(*)$.  
\item  \emph{We show that $f_3 \not \in \Sigma_2$.} \\
If $f_3\in \Sigma_2$, the representation $\rho$ is trivial, hence a contradiction. Therefore,  $f_3 \not \in \Sigma_2$ and one can define an integer $a \leq m_3$ such that $X^a$ is the monomial of $f_3$ with highest degree among those that do not belong to $\Sigma_2$. Since $f_3$ is assumed to be reduced mod $\wp(k[X])$, then $a \neq 0$ mod $p$.
\item \emph{We show that $p$ divides $a-1$.} \\
 Consider the equation:
\begin{equation} \label{f3}
\forall \, y \in V, \quad \Delta_y(f_3)=\ell_{1,3}(y) \, f_1(X)+ \ell_{2,3}(y) \, f_2(X)  \quad \mod \wp(k[X])
\end{equation}
where $\ell_{1,3}$ and $\ell_{2,3}$ are non zero linear forms from $V$ to $\F_p$.
The monomials of $f_3$ with degree strictly lower than $a$ belong to $\Sigma_2$. So they give linear contributions in $\Delta_y(f_3)$ mod $\wp(k[X])$ (cf. \cite{MR2} Lemma 3.9). Assume that $p$ does not divide $a-1$. Then, for all $y$ in $V$, equation \eqref{f3} gives the following equality mod $\wp(k[X])$:
$$c_a(f_3) \, a \, X^{a-1} +\mbox{lower degree terms}=(\ell_{1,3}(y)+\gamma_2 \, \ell_{2,3}(y))\, X^{1+p^s}  +\mbox{lower degree terms} $$
where $c_a(f_3) \neq 0$ denotes the coefficient of $X^a$ in $f_3$.
If $a-1 >1+p^s$, then $y=0$ for all $y$ in $V$ and $V=\{0\}$ which is excluded for a big action (cf. \cite{MR} Prop. 2.2).
If $a-1 <1+p^s$, then, $\ell_{1,3}(y)+\gamma_2 \, \ell_{2,3}(y)=0$, for all $y$ in $V$. It follows that $\gamma_2 \in \F_p$, which is another contradiction. So, $a-1=1+p^s$ and by equating the corresponding coefficients in \eqref{f3}, one gets: $a\, y= \ell_{1,3}(y)+\gamma_2 \, \ell_{2,3}(y)$, for all $y$ in $V$. So, $V \subset \F_p+\gamma_2 \, \F_p$ and $v \leq 2$. As $v=2\,s$, we deduce that $s=1$, which is a contradiction.
Thus, $p$ divides $a-1$ and one can write $a=1+\lambda \, p^t$ with $t\geq 1$ and $\lambda \geq 2$, because of the definition of $a$. We also define $j_0:= a-p^t$. 
\item \emph{We show that $v \geq t+1$.} \\
By \cite{MR2} (Lemma 3.11), $p^v \geq m_3+1 >m_3-1 \geq a-1=\lambda\, p^t \geq 2\, p^t$. 
This implies $v \geq t+1$.
\item \emph{We show that $j_0=1+p^s$.} \\
If $j_0<1+p^s$, we gather the same contradiction as the one found in \cite{MR2} [proof of Theorem 5.6, point 4, with $i=2$]. Now, assume that $j_0 >1+p^s$. As in \cite{MR} [proof of Theorem 5.1, point 6], we prove that the coefficient of $X^{j_0}$ in the left-hand side of \eqref{f3} is $T(y)$, where $T$ is a polynomial of $k[X]$ with degree $p^t$. If $j_0 >1+p^s$, then $T(y)=0$, for all $y$ in $V$. This implies $V \subset Z(T)$ and $v \leq t$, which contradicts the previous point.
\item \emph{We show that either $v=t+1$ or $v=t+2$.} \\
We have already seen that $v \geq t$. As $j_0=1+p^s$, we equate the corresponding coefficients in \eqref{f3} and obtain $T(y)= \ell_{1,3}(y)+\gamma_2 \, \ell_{2,3}(y)$, for all $y$ in $V$. As $\ell_{1,3}(y) \in \F_p$ and $\ell_{2,3}(y) \in \F_p$, we get $T(y)^p-T(y)= \ell_{2,3}(y) \, (\gamma_2^p-\gamma_2)$, with $\gamma_2 \not \in \F_p$. Then, for all $y$ in $V$, $R(y):= \frac{T(y)^p-T(y)}{\gamma_2^p-\gamma_2} =\ell_{2,3} (y) \in \F_p$ and $V \subset Z(R^p-R)$. In particular, $v \leq t+2$.
\item \emph{We show that $m_3=a=1+p^s+p^t$.} \\
Assume that $m_3 >a$. Then, by definition of $a$, $m_3=1+p^{s_3}$ with $s_3 \geq s$. Note that $s_3 \geq s+1$.
Otherwise, $m_3=1+p^s =j_0 <a$. On the one hand, $|G|=p^{3+v}=p^{3+2s}$. On the other hand,
 $$g =\frac{p-1}{2} \, (p^s+p^{s+1}+p^2(m_3-1)) =  \frac{p-1}{2} \,p^s \, (1+p+p^{2+s_3-s}) \geq \frac{p-1}{2} \, p^s \, (1+p+p^3)$$
Thus, $\frac{|G|}{g^2} \leq \frac{4}{(p^2-1)^2} \, \frac{p^3\, (p+1)^2}{(1+p+p^3)^2} <\frac{4}{(p^2-1)^2}$.
This contradicts condition $(*)$, so $m_3=a$.
\item \emph{We show that $s=2$ and $v=4$. In particular, $\gamma_2 \in \F_{p^2}-\F_p$.} \\
We already know that $s\geq 2$ and $v=2\,s \geq 4$. So, $|G| =p^{3+v} \leq p^7$. Assume that $s \geq 3$. Then, as $t \geq 1$, we get:
$g =\frac{p-1}{2} \, (p^s+p^{s+1}+p^2(m_3-1))=\frac{p-1}{2} \, (p^s+p^{s+1}+p^{s+2}+p^{t+2}) \geq 
\frac{p-1}{2} \, (2\,p^{3}+p^4+p^5)$. It follows that 
$\frac{|G|}{g^2} \leq \frac{4}{(p^2-1)^2} \, \frac{p\, (p+1)^2}{(2+p+p^2)^2} <\frac{4}{(p^2-1)^2}$, which is a contradiction. So $s=2$ and $v=4$. We have previously mentionned that $\gamma_2 \in \F_{p^d}-\F_p$, where $d$ is an integer dividing $s$. As $s=2$, the only possibility is $d=2$.
\item \emph{We deduce that $t=s=2$, so $m_3=1+2\, p^2$ and $p \geq 5$. } \\
We have seen $v=t+1$ or $v=t+2$, with $t \geq 1$. As $v=4$, there are two possibilities either $t=2$ or $t=3$.
If $t=3$, $|G|=p^7$ and $g=\frac{p-1}{2} \,p^2\,  (1+p+2\, p^{3})$. So, 
$\frac{|G|}{g^2} \leq \frac{4}{(p^2-1)^2} \, \frac{p^3\, (p+1)^2}{(1+p+2\,p^3)^2} <\frac{4}{(p^2-1)^2}$.
Therefore, $t=2=s$. In this case, $\frac{|G|}{g^2} = \frac{4}{(p^2-1)^2} \, \frac{p^3\, (p+1)^2}{(1+p+2\,p^2)^2}$ and condition $(*)$ requires $p \geq 5$.
\item \emph{We gather the parametrization of $f_1$, $f_2$ and $V$. } \\
As $s=d=2$, $f_1$ reads $f_1(X)=X\, S_1(X)$ with $S_1(F)= \sum_{j}^{s/d} \, a_{jd} \, F^{jd}= a_0 \, I+F^2$, since $S_1$ is assumed to be monic. Then, 
$$f_1(X)=X\, (X^{p^2}+a_2\, X^2) \quad \mbox {and} \quad f_2(X)=\gamma_2 \, X\, (X^{p^2}+a_2\, X^2)+b_1\, X$$
with $a_2\in k$, $b_1\in k$ and   $\gamma_2 \in \F_{p^2}-\F_p$.
In this case, $$V=Z(Ad_{f_1})= Z(X^{p^4}+ + 2\, a_2^{p^2}\, X^{p^2}+X)$$
\item \emph{We show that $f_3 \in \Sigma_4$ but $f_3 \not \in \Sigma_4-\Sigma_3$. } \\
By \cite{MR2} (Thm. 3.13), $f_3 \in \Sigma_4$. We now show that $f_3$ does not have any monomial in $\Sigma_4-\Sigma_3$. Indeed, as $m_3=1+2\, p^2$, the possible monomials of $f_3$ in $\Sigma_4-\Sigma_3$ are $X^{1+2\,p+p^2}$, $X^{2+p+p^2}$, $X^{3+p^2}$, $X^{1+3\,p}$,
$X^{2+2\,p}$, $X^{3+p}$ and $X^4$. Now, equate the coefficients of the monomial $X^{1+p+p^2} \in \Sigma_3$ in each side of \eqref{f3}. In the left-hand side, i.e. in $\Delta_y(f_3)$ mod $\wp(k[X])$, $X^{1+p+p^2}$ is produced by monomials $X^b$ of $f_3$ that belong to $\Sigma_4-\Sigma_3$ and satisfy $b> 1+p+p^2$. This leaves only two possibilities:
$X^{1+2\,p+p^2}$ and $X^{2+p+p^2}$. In the right-hand side of \eqref{f3}, $X^{1+p+p^2} \in \Sigma_3-\Sigma_2$ does not occur since $\ell_{1,3}(y) \, f_1(X)+ \ell_{2,3}(y) \, f_2(X)$ lies in $\Sigma_2$. It follows that, for all $y$ in $V$, 
$2\,  c_{2+p+p^2}  \, y^p+2\,  c_{1+2\, p+p^2}\, y=0$, where $c_t$ denotes the coefficient of the monomial $X^t$ in $f_3$. As $v=dim_{\F_p} \, V=4$, we deduce that
$c_{2+p+p^2} =c_{1+2\, p+p^2}=0$. 
We go on this way and equate successively the coefficients of $X^{2+p^2}$, $X^{1+2\,p}$, $X^{2+p}$ and $X^3$ to prove that $f_3$ does not contain any monomial in $\Sigma_4-\Sigma_3$. Therefore, $f_3$ reads as follows:
 $$f_3(X)= c_{1+2\,p^2} \, X^{1+2\,p^2} + c_{1+p+p^2} \, X^{1+p+p^2} + c_{2+p^2} \, X ^{2+p^2} + c_{1+p^2} \, X ^{1+p^2} + $$
$$c_{1+2\,p^2} \, X^{1+2\,p} + c_{2+p} \, X^{2+p}+
c_{1+p} \, X^{1+p} + c_{3} \, X^3+  c_{2} X^2 + c_{1} \,  X$$
 \item \emph{We determine $f_3$. } \\
We finally have to solve \eqref{f3} with $f_1$, $f_2$ and $f_3$ as described above.
Calculation show that $c_{1+p+p^2}=c_{1+2\,p^2}= c_{2+p}=0$ and that the coefficients $a_2$, $c_{2+p^2}$ and $c_{3}$ can be expressed as rational functions in $c_{1+2\,p^2}$ (see formulas in the table given in the proposition). To conclude, one has to distinguish the cases $b_1 \neq 0$ and $b_1=0$.
In the first case, $b_1$, $c_{1+p}$ and $c_{2}$ can be expressed as rational functions in $c_{1+p^2}$ whereas $e:=c_{1+p^2}-c_{1+p^2}^{p^2}$ belongs to the set of zeroes of a polynomial whose coefficients are rational functions in $c_{1+2\,p^2}$ (see table). When $b_1=0$, then $c_{1+p}=0$, $c_{1}= c_{1+p^2} \, a_2$ and $c_{1+p^2} \in \F_{p^2}$.
It follows that $f_3(X)= c_{1+2\,p^2}\,  X^{1+2\,p^2} + c_{2+p^2} \,  X ^{2+p^2} + c_{3}\,  X^3+ \, c_{1} \,  X + c_{1+p^2} \, f_1(X)$. As $\gamma_2\in \F_{p^2}-\F_p$,  $\{1, \gamma_2\}$ is a basis of $\F_{p^2}$ over $\F_p$. Write $\epsilon =\epsilon_1+\epsilon_2\, \gamma_2$, with $\epsilon_1$ and $\epsilon_2$ in $\F_p$. By replacing $f_3$ with $f_3-(\epsilon_1\, f_1+\epsilon_2\, f_2)$, one obtains the expected formula. 
\end{enumerate}
\item The equality $G=A_{\infty,1}$ derives from the maximality of $V=Z(Ad_{f_1})$ (see Proposition \ref{prop11}). The unicity of the $p$-Sylow subgroup is due to Remark \ref{italien}. $\square$
\end{enumerate}
\medskip

\noindent The last case: $\ell_{1,2}\neq 0$ and $\ell_{2,3}\neq 0$, generalizes the results obtained in \cite{MR2} (section 6.2).
\medskip

\begin{proposition} \label{p3cas3}
Let $(C,G)$ be a big action satisfying $\mathcal{G}_*^{p^3}$ such that $[G',G] \neq \{e\}$. 
We keep the notations introduced in section 5.1. Assume that $\ell_{1,2}\neq 0$ and $\ell_{2,3}\neq 0$. 
\begin{enumerate}
\item Then, $p \geq 11$ and there exists a coordinate $X$ for the projective line $C/G_2$ such that
the functions $f_i$'s can be parametrized as follows:
\medskip

\noindent \begin{tabular}{|c|c|}
\hline
 $f_1$  &  $f_1(X)=X^{1+p}+a_2\, X^2$ \\
\hline 
$V$ & $V= Z(Ad_{f_1}) =Z(X^{p^2}+2\, a_2^p \,X^p +X)$\\
\hline
\hline
$f_2 $ & $f_2(X)= b_{1+2\,p} \, X^{1+2\,p} + b_{2+p} \, X^{2+p}  + b_{3} \,X^3+  b_{1}\,X$\\
\hline
$b_{1+2\,p}$ & $ b_{1+2\,p} \in k^{\times}$ \\
\hline
$a_2$ & $2\, a_2^p= -b_{1+2\,p}^{-p} \, (b_{1+2\,p}+b_{1+2\,p}^{p^2}) \Leftrightarrow b_{1+2\,p} \in V$\\
\hline
$V$ &$ V=Z(X^{p^2}-b_{1+2\,p}^{-p} \, (b_{1+2\,p}+b_{1+2\,p}^{p^2} ) \,X^p +X)$\\
\hline
$b_{2+p}$& $ b_{2+p}=-b_{1+2\,p}^p$\\
\hline
$b_{3}$ & $  3 \,b_{3}^p= b_{1+2\,p}^{-p} \, (b_{1+2\,p}^{2\,p^2}-b_{1+2\,p}^{2})$ \\
\hline
$\ell_{1,2}$ & $\ell_{1,2}(y)=2\, (b_{1+2\,p} \, y^p-b_{1+2\,p}^p \, y)$\\
\hline
\hline
$f_3 $ & $f_3(X)= c_{1+3\,p} \, X^{1+3\,p} + c_{2+2\,p} \, X^{2+2\,p}  + c_{1+2\,p} \, X^{1+2\,p} +c_{3+p} \, X^{3+p} $\\
& $+c_{2+p} \, X^{2+p} +c_{1+p}\, X^{1+p} +c_{4} \, X^4+ c_{3} \,X^3+ c_{2} \, X^2+ c_1\,X$\\
\hline
$c_{1+3\,p}$ & $ 3\, c_{1+3\,p}=2\, b_{1+2\,p}^2 $ \\
\hline
$c_{2+2\,p}$ & $c_{2+2\,p}=-b_{1+2\,p}^{1+p}$ \\
\hline
$c_{3+p}$& $ 3\, c_{3+p}=2 \,b_{1+2\,p}^{2\,p}$\\
\hline
$c_{4}$ & $6\,c_{4}^p=-b_{1+2\,p}^{-p} \, (b_{1+2\,p}^3+b_{1+2\,p}^{3\,p^2})$ \\
\hline
$c_{1+2\,p} $ & $ c_{1+2\,p} \in V $\\
\hline
$c_{2+p}$ & $c_{2+p}=-c_{1+2\,p} ^p$ \\
\hline
$c_{3}$ & $3\,c_{3}^p=b_{1+2\,p}^{-p} \, (b_{1+2\,p}+b_{1+2\,p}^{p^2})\, (c_{1+2\,p} ^{p^2}-c_{1+2\,p})$\\
\hline
$c_{1+p}$ & $c_{1+p} \in k$\\
\hline
$b_{1} $ & $2\, b_{1}^p = b_{1+2\,p}^{-p} \, ( c_{1+p}^p -c_{1+p}) $\\
\hline
$c_{2}$ & $2\, c_{2}^p= -b_{1+2\,p}^{-p} \, (c_{1+p}^p \, b_{1+2\,p}^{p^2} +c_{1+p} \, b_{1+2\,p})$ \\
\hline
$c_1$ & $c_1 \in k$\\
\hline
$\ell_{1,3}$&   $\ell_{1,3}(y)=2\,( c_{1+2\,p} \, y^p\, c_{1+p}^p -y)+ 2\, b_{1+2\,p}^2 \, y^{2\,p}-4\, b_{1+2\,p}^{1+p} \, y^{1+p} +2\, b_{1+2\,p}^{2\,p} \, y^2 $\\
&$=2\,( c_{1+2\,p} \, y^p-c_{1+p}^p \, y) +\ell_{1,2}^2(y)/2$ \\
\hline
$\ell_{2,3}$ & $\ell_{2,3}(y)=2\, (b_{1+2\,p} \, y^p-b_{1+2\,p}^p \, y)$\\
\hline
\end{tabular} 

\bigskip

\noindent Therefore, the solutions can be parametrized by 3 algebraically independent variables over $\F_p$, namely $b_{1+2\,p}\in k^{\times}$, $c_{1+p} \in k$ and $c_1 \in k$. One also finds a fourth parameter $c_{1+2\,p}$ which runs over $V$. So, for a given $b_{1+2\,p}$, the parameter $c_{1+2\,p}$ takes a finite number of values. \\ 
\medskip
\noindent Moreover,
$$\frac{|G|}{g}=\frac{2\,p}{p-1} \, \, \frac{p^3}{1+2\,p+3\, p^2} \quad \mbox{and} \quad \frac{|G|}{g^2}=\frac{4}{(p^2-1)^2} \,\,  \frac{p^3\, (p+1)^2}{(1+2\,p+3\,p^2)^2}$$ 
\item  $G=A_{\infty,1}$ is the unique $p$-Sylow subgroup of $A$. Furthermore, $Z(G)$ is cyclic of order $p$. 
\end{enumerate}
\end{proposition}

\noindent \textbf{Proof:}
\begin{enumerate}
\item 
In this case, the group $G$ satisfies the third condition of \cite{MR2} (Prop. 5.2). So, we deduce from
\cite{MR2} (Thm. 5.6) that $m_1=1+p^s$, $m_2=1+2\, p^s$, $m_3=1+3\, p^s$ with $p\geq 5$ and $v=s+1$.
Furthermore, it follows from Lemma \ref{G12} that $s \in \{1,2\}$. Assume that $s=2$. Then, $|G|=p^6$,
$g=\frac{p-1}{2} \, p^2\, (1+2\,p+3\, p^2)$, so 
$\frac{|G|}{g^2}=\frac{4}{(p^2-1)^2} \, \frac{p^2\, (p+1)^2}{(1+2\,p+3\,p^2)^2} <\frac{4}{(p^2-1)^2}$.
This is a contradiction, hence $s=1$. In this case, 
$\frac{|G|}{g^2}=\frac{4}{(p^2-1)^2} \, \frac{p^3\, (p+1)^2}{(1+2\,p+3\,p^2)^2}$ and condition $(*)$ is satisfied as soon as $p\geq 11$.
Then, we deduce from Proposition \ref{casp2nontrivial} the parametrization of $f_1$, $V$ and $f_2$ mentionned in the table.
Besides, we deduce from \cite{MR2} (Thm. 5.6) that $f_3$ is in $\Sigma_4-\Sigma_3$ with $m_3=1+3\,p$. This means that $f_3$ reads as follows:
$$f_3(X)= c_{1+3\,p} \, X^{1+3\,p} + c_{2+2\,p} \, X^{2+2\,p}+ c_{1+2\,p} \, X^{1+2\,p}+
c_{3+p} \, X^{3+p}$$
$$+ c_{2+p} \, X^{2+p}+c_{1+p} \, X^{1+p} + c_{4} \, X^4 
+ c_{3} \, X^3 + c_{2} \, X^2 + c_1 \, X $$
We determine the expressions of the coefficient by solving the equation:
$$ \forall \, y \in V, \quad \Delta_y(f_3)= \ell_{1,3}(y) \, f_1(X)+ \ell_{2,3} \, f_2(X) \quad
\mod \wp(k[X])$$
with $\ell_{1,2}(y)=\ell_{2,3}(y)= 2\, (b_{1+2\,p} \, y^p- b_{1+2\,p}^p \, y)$ (cf. \cite{MR2}, Prop. 5.4.1). The results are gathered in the table above. 
\item The equality $G=A_{\infty,1}$ derives from \cite{MR} (Cor. 5.7). The unicity of the $p$-Sylow subgroup comes from  Remark \ref{italien}. The description of the center is due to \cite{MR} (Prop. 6.15). $\square$
\end{enumerate}

\bigskip

\noindent {\bf \Large{Acknowledgements}.}\\
I wish to express my gratitude to my supervisor, M. Matignon, for his valuable advice, his precious help and his permanent support all along the elaboration of this paper.

\end{document}